\documentclass[11pt] {article} 
  \usepackage{amsmath} 
  \usepackage{mathrsfs}
  \usepackage[backref]{hyperref}

\newtheorem{theorem}{Theorem}[section]
\newtheorem{lemma}{Lemma}[section]

\newcommand{\eqnsection}{
   \renewcommand{\theequation}{\thesection.\arabic{equation}}
   \makeatletter
   \csname @addtoreset\endcsname{equation}{section}
   \makeatother}
 
\def \be{\begin{equation}}
\def \ee{\end{equation}}
\def \bt{\begin{theorem}} 
\def \et{\end{theorem}}
\def \bl{\begin{lemma}} 
\def \el{\end{lemma}} 
\def \bea{\begin{eqnarray}}
\def \eea{\end{eqnarray}}
\def \bas{\begin{eqnarray*}}
\def \eas{\end{eqnarray*}}

\def \al{\alpha} 
\def \bb{\beta}
\def \ga{\gamma} 

\def \de{\delta}
\def \De{\Delta} 
\def \ep{\epsilon}

\def \la{\lambda} 

\def \ka{\kappa}
\def \om{\omega}

\def \si{\sigma}

\def \ff{\infty}
\def \wh{\widehat}
\def \wt{\widetilde}

\def \rar{\rightarrow}

\def \CC{{\cal C}}

\def \EE{{\cal E}}
\def \FF{{\cal F}}
\def \GG{{\cal G}}
\def \HH{{\cal H}}

\def \JJ{{\cal J}}
\def \KK{{\cal K}}
\def \LL{{\cal L}}
\def \MM{{\cal M}}

\def \PP{{\cal P}}

\def \VV{{\cal V}}

\def \({\left(}
\def \){\right)}

\def \lc{\left\{}
\def \rc{\right\}}

\def \nn{\nonumber}

\def \bc{\begin{center} }
\def \ec{\end{center} }
\def \bs{\begin{slide} }
\def \es{\end{slide} }

\def\square{{\vcenter{\vbox{\hrule height.3pt
        \hbox{\vrule width.3pt height5pt \kern5pt
           \vrule width.3pt}
        \hrule height.3pt}}}}
\def\qed{{\hfill $\square$ \bigskip}}

\def \bpr{\begin{proof} }
\def \epr{\end{proof} }

\def \st{\stackrel{def}{=}}

\eqnsection
\begin{document}

\title{Intersection local times for interlacements}

    \author{ Jay Rosen \thanks{Research  was  supported by  grants from the National Science Foundation.}}
\maketitle

\begin{abstract}   
We define renormalized intersection local times for random  interlacements of
 L\'evy processes  in $R^{d}$ and prove an isomorphism theorem relating   renormalized intersection local times with associated Wick polynomials.  
 \end{abstract}

\bibliographystyle{amsplain}

  \section{Introduction}   \label{sec-1}

  A random interlacement  in $R^{d}$ is a particular Poison process $\mathcal{I}_{\al}$ of paths  in $R^{d}$, \cite{Sz1,Sz2, Sz3}. We
 consider the  $n$-fold intersections of random interlacements of
 L\'evy processes in $R^{d}$. This entails studying  functionals of the form
\be
 \ga_{n,\ep}(\nu)
\stackrel{def}{=}\int  \(\sum_{\om\in\mathcal{I}_{\al}} \int f_{\ep}(Y_{t}(\om)- x)\,dt \)^{n}
\,d\nu(x),\label{3}
\ee  
where $Y_{t}(\om)=\om(t)$ for a path $\om\subset R^{d} $, $f_{\ep}$ is an approximate $\delta$-function at zero and $\nu$ is a finite  measure on $R^{d}$.
 Ideally we would like to take the limit   of $ \ga_{n,\ep}(\nu)$ as $\ep$ goes to 0, but  in general   the limit is
infinite   for all  $n\ge 2$.
To  deal with this we use a technique  called renormalization, which consists of
forming a
linear combination of  the 
$\{\ga_{k,\ep}(\nu)\}_{k=1}^n$ which has a finite limit, $ L_{n}(\nu)$, as $\ep\to 0$.   We study  the behavior of $L_{n}(\nu)$  as a function of $\nu$.

Renormalized  intersection local time (rilt) for Markov processes originated with the work of  Varadhan
\cite{Va} who studied planar Brownian  rilt for its role in quantum field theory. 
Renormalized intersection local time turns out to be the right tool  for the
solution of certain ``classical'' problems such as the asymptotic expansion of 
the area of the Wiener sausage in the plane and the range of
    random walks, \cite{BR}, \cite{LGa}, \cite{LGb}. For further work on
rilt see Dynkin
\cite{Da},   Bass and Khoshnevisan
\cite{BKb}, Rosen
\cite{Rb} and Marcus and Rosen \cite{MRmem}.

We set 
$L_{1,\ep}(\nu)=\ga_{1,\ep}(\nu)$ and define   recursively  
\begin{equation}
L_{n,\ep}(\nu)= \ga_{n,\ep}(\nu)-\sum_{j=1}^{n-1}c_{n,j,\ep}L_{j,\ep}(\nu)\label{zzx}
\end{equation}
where the $c_{n,j,\ep}$ are constants which diverge as $\ep\rar 0$; see (\ref{rilt.1})   and (\ref{rilt.2b}).  
We show that for a wide class of random interlacements and   finite  compactly supported measures $\nu$ on $R^{d}$ 
\begin{equation}
L_{n}(\nu):=\lim_{\ep \rar 0}L_{n,\ep}(\nu)\quad\mbox{    exists in all  $L^{p}$.}\label{rilt.9intro}
\end{equation}
We refer to $L_{n}(\nu)$ as the $n$-fold renormalized  intersection
local time   of $\mathcal{I}_{\al}$ with respect to   $\nu $. 

As indicated, a random interlacement $\mathcal{I}_{\al}$ is a Poisson process of paths associated with a transient Markov pocess. A Poisson process is determined by its intensity measure. For random interlacements the intensity measure $\al\mu_{m}$, $\al>0$, is a measure on bilateral paths which measures  geometric properties of the paths, rather than their particular parametrization. More precisely,  $\mu_{m}$ is invariant under time shifts. Before giving the precise characterization of $\mu_{m}$ and providing references for further details, let us give some indication of how $\mu_{m}$ looks for Brownian motion in $R^{3}$. Let $K$ be a compact subset of $R^{3}$. 
Consider the set $\mathcal{A}$ of  paths $X$ which, up to time shift, hit $K$ for the first time at $t=0$, lie in some set $\mathcal{A}^{+}\subseteq C(R_{+}^{1},R^{3} )$ for $t\geq 0$, and 
 $\mathcal{A}^{-}\subseteq C(R_{-}^{1},R^{3} )$ for $t\leq 0$. Then we want
 \begin{equation}
\mu_{m}(\mathcal{A})=\int P^{x}(\mathcal{A}^{+})P_{K}^{x}(\mathcal{A}^{-}\circ r)\, e_{K}(dy). \label{exam.1}
 \end{equation}
 Here $e_{K}(dy)$ is the equilibrium measure of $K$ for Brownian motion in $R^{3}$, $P^{x}$ is the usual  probability for (one-sided)  Brownian paths starting at $x$, $P_{K}^{x}$ is $P^{x}$ conditioned never to return to $K$ and $r(\om)(t)=\om(-t)$. What follows is a more precise characterization  of $\mu_{m}$.

In this paper we deal only with random interlacements  of symmetric  L\'evy processes in $R^{d}$. However,  random interlacements can be defined quite generally.  Let 
$X=\left(\Omega,\mathcal{F},\mathcal{F}_t,X_t, \theta_t,P^x
\right)$ be a transient Borel right process with a locally compact state space $S$ and  
potential
densities $u(x,y)$ with respect to a $\si$-finite excessive measure $m$. We use $P_{t}$ to denote the semigroup for $X$. We assume that $m$ is dissipative, that is, that $\int u(x,y)f(y)\,m(dy)<\ff$ $m$-a.e. for each non-negative  $f\in L^{1}(m)$.

Let $W$ denote the set of paths $\om: R^{1}\mapsto S\cup \De$ which are $S$ valued and right continuous on some open interval $(\al (\om), \beta (\om))$ and $\om (t)=\De $ otherwise. Let $Y_{t}=\om (t)$,   and define the shift operators
\begin{equation}
(\si_{t}\om) (s)=\om (t+s),\hspace{.2 in}s,t\in R^{1}.\label{q.1}
\end{equation}
Set $\mathcal{F}=\si\(Y_{s}, s\in R^{1}\)$ and $\mathcal{F}_{t}=\si\(Y_{s}, s\leq t\)$.
Let $\mathcal{A}$ denote the $\si$-algebra of shift invariant events in $\mathcal{F}$. 
 The quasi-process associated with $X$ is the measure $\mu_{m}$ on $\(W,\mathcal{A} \)$ which satisfies the following two conditions: 
 \begin{equation}
(i): \hspace{.5 in}\mu_{m}\(\int_{R^{1}}f\(Y_{t}\)\,dt\)=m(f)\label{q.2}
 \end{equation}
 and (ii): if  $T$ is any intrinsic stopping time, then $Y_{T+t}, t>0$ is Markovian with semigroup $P_{t}$ under  
 $\mu_{m}|_{\{T\in R\}}$. 
 An $\mathcal{F}_{t^{+}}$ stopping time $T$ is called intrinsic if $\al\leq T\leq\bb$ on $\{T<\ff\}$ and 
 $T=t+T\circ \si_{t}$ for all $t\in R^{1}$. Note that since our times run over $R^{1}$, this definition does not contain the usual condition that $T\geq t$. If $T<t$, then   $T\circ \si_{t}$ will be negative.  A first hitting time is an example of an intrinsic stopping time.  The quasi-process associated with $X$ will exist under the conditions of the previous paragraph,   see \cite[XIX]{DM4}. The name `quasi-process' refers to the fact that $\mu_{m}$ is only defined on the $\si$-algebra $\mathcal{A}$ of shift invariant sets, hence for example,  if  $B\subseteq S$, one cannot ask for the measure of the set $\{   Y_{t}\in B\}$.
 
Random interlacements are the `soup' of a quasi-process. More precisely, for any $\al>0$, the random interlacement $\mathcal{I}_{\al}$ associated with $X$ is the Poisson process in $W$ with intensity measure $\al\mu_{m}$. We let 
$\mathcal{P}_{\al}$ denote probabilities for the process $\mathcal{I}_{\al}$.

In the rest of this paper we take $X=\{X(t),t\in R^+\}$  to be a symmetric   L\'evy process in $R^d$  that is either transient to begin with or made transient by killing at the end of an  independent exponential time with mean $1/\ka$, and we take Lebesgue measure $dx$ as our $\si$-finite excessive measure $m$. For simplicity, we assume that $X$ is radially symmetric with characteristic exponent of the form $\psi(|\cdot |)$, i.e.
\begin{equation}
E\(e^{i\xi X_{t}}\)=e^{-t \psi(|\xi|)},\label{sl.00}
\end{equation}  
 where $\psi$   is regularly varying at infinity with index $\bb$ satisfying
 \be  \(1-\frac{1}{2n}\)d<\bb\le d.\label{8.34} 
 \ee
We assume further that  $1/\psi(|\cdot |)$ is locally integrable, (this is automatic if we are dealing with an exponentially killed process), and
    \begin{equation}
      \int_{R^{d}}  1/\psi(|\xi|)  \,d\xi=\ff,\label{sl1.5}
   \end{equation}
which is automatic if $\bb<d$.
We mention that Lebesgue measure $dx$ is dissipative for such processes, hence the random interlacement $\mathcal{I}_{\al}$ associated with such  $X$ exist.

 The potential densities $u(x,y)=u(x-y)$,  which may be infinite on the diagonal,  are always weakly positive definite.
We use $\GG^{m}$ to denote the set of  positive measures $\nu$   which are bounded with bounded potential and  for which
\be
\int \int (u(x,y))^{m}\,d\nu(x)\,d\nu(y)<\ff.\label{14.2} 
\ee 
For compact $K\subseteq R^{d}$ 
we use $\GG_{K}^{m}$ to denote the set of $\nu\in \GG^{m}$ with support in $K$.
It should be
understood that when we say $\nu\in\GG_{K}^{m}$, that this is with respect to the
potential of some given L\'evy process. We point out that for all $n\geq 2$, $\GG_{K}^{2n}$ will be empty unless $d=1$ or $2$.

 \bt\label{theo-multiriltintro-m} Let $\mathcal{I}_{\al}$ be the random interlacement   associated with
 a L\'evy process  as above.
 Let $\nu$ be a finite measure in $\mathcal{G}_{K}^{2n}$, for some compact $K\subset  R^d$.
Then (\ref{rilt.9intro}) holds.
   \et

Our main result is an isomorphism theorem relating the renormalized intersection local times $L_{n}(\nu)$ with associated Wick polynomials. The definition of Wick powers $:G^{j}:(\nu)$ is recalled in Section \ref{sec-prelim} and that of the `mixed terms' $\({ :G^{2j}:\over 2^{j}}\times L_{n-j}\)(\nu)$
is given in (\ref{cross1}).
  
    \bt[Isomorphism Theorem]\label{theo-ILT}
For any  $\al >0$, compact $K\subset S$ and countable  set  $D\subseteq  \mathcal{G}_{K}^{2n}$,
\bea
&&
\Big\{\sum_{j=0}^{n}{ n\choose j}\({ :G^{2j}:\over 2^{j}}\times L_{n-j}\)(\nu),\,\nu\in D, \mathcal{P}_{\al^{2}}\times P_{G}\Big\}\label{1.0}\\
&&\stackrel{law}{=}
\Big\{\sum_{j=0}^{2n}{ 2n\choose j}\al^{(2n-j)}{ :G^{j}:(\nu)\over 2^{j/2}},\,\nu\in D,P_{G} \Big\}.\nn
\eea
\et

 In particular, when $n=2$ this says that
for any  $\al >0$, compact $K\subset S$ and countable set $D\subseteq  \mathcal{G}_{K}^{4}$,
\bea
&&
\Big\{ L_{2}(\nu) +2\({ :G^{2}:\over 2}\times L_{1}\)(\nu)+{ :G^{4}:(\nu)\over 2^{2}},\,\nu\in D, \mathcal{P}_{\al^{2}}\times P_{G}\Big\}\nn\\
&&\stackrel{law}{=}
\Big\{{ :G^{4}:(\nu)\over 2^{2}}+4\al  { :G^{3}:(\nu)\over 2^{3/2}}+6\al^{2}  { :G^{2}:(\nu)\over 2} \nn\\
&& \hspace{1.5 in}+4\al^{3} { G(\nu)\over 2^{1/2}}+ \al^{4}  |\nu|,\,\nu\in D,P_{G} \Big\}.\label{2.0}
\eea
 
 At first glance these isomorphism theorems may seem too complicated to work with. However, we have found similar isomorphism theorems very useful, see \cite{R} and especially 
  \cite[p. 33]{MRmem}.
  Here is a particularly straightforward example which mirrors our results in \cite{MRmem} for ordinary intersection local times of L\'evy processes. 
    We are concerned with the continuity of  $\{L_{n}(\nu),\nu\in \VV\}$, where $\VV$ is some metric space.
    
\medskip	   Let
$\tau_{2n}(\xi)$ denote the Fourier transform of $(u(x))^{2n}$ so that 
\be
\int\tau_{2n}(\xi)|\hat\nu(\xi)|^2\,d\xi=\int\!\! \int \(u(x,y)\)^{2n}
\,d\nu(x)\,d\nu(y) .\label{a5.25}
\ee    For any finite positive measure $\nu$ on
$R^d$, let $\nu_{x}(A)=\nu(A-x)$.

\bt\label{theo-cont14.1b}  Under the hypotheses of Theorem  \ref{theo-multiriltintro-m}, if 
$\nu\in \GG_{K}^{2n}$ is such that
\be
\int_{1}^{\ff}  \(\int_{|\xi|\geq x}\tau_{2n}(\xi)|\hat{\nu}(\xi)|^{2}\,d\xi\) ^{1/2} \frac{
 (\log x)^{n-1}}{x}\,dx<\ff,\label{b101}
\ee   
then 
$\{L_{n}(\nu_x),\,x\in R^m\}$ is continuous almost surely.
\et

In particular, for the random interlacement   associated with exponentially killed 
Brownian motion in
$R^{2}$, this is the case when \be
|\hat\nu(\xi)|=O\(\frac1{(\log|\xi|)^{2n+\ep}}\)\qquad as\quad  |\xi|\to\ff.\label{xx}
\ee

Furthermore for the random interlacement   associated with a L\'evy process $X$ in $R^2$  with L\'evy exponent
asymptotic to
$\la^2/(\log |\la|)^a$, $a>0$, as $\la\to\ff$, (see \cite[p. 5]{MRmem}), 
$\{L_{n}(\nu_x),\,x\in R^m\}$  is continuous almost surely if (\ref{xx}) holds
with $2n$ replaced by
$2n(1+a/2)$.

The $n=1$ case of our theorem,  
\bea
&&
\Big\{\textstyle{ 1\over 2}:G^{2}:(\nu)+ L_{1}(\nu),\,\nu\in D, \mathcal{P}_{\al^{2}}\times P_{G}\Big\}\label{0.0jh}\\
&&\stackrel{law}{=}
\Big\{\textstyle{ 1\over 2}:G^{2}:(\nu)+\sqrt{2}\al G(\nu)+\al^{2}  |\nu|,\,\nu\in D,P_{G} \Big\},\nn
\eea
is essentially due to Sznitman. Formally, it says that there is an equivalence in law  between
\begin{equation}
G^{2}/2 +L_{1}\hspace{ .2in}\mbox{and}\hspace{ .2in}\( G/2^{1/2}  +\al\)^{2}.\label{0.1jh}
\end{equation}
Taking the $n$'th power of both sides suggest there should be an equivalence in law  between
\begin{equation}
 \({ G^{2}\over 2} +L_{1} \)^{n}=\sum_{j=0}^{n} \binom {n}{ j} \({ G^{2j}\over 2^{j}}\times L_1^{n-j}\)\label{0.2jh}
\end{equation}
and
\begin{equation}
 \({ G\over 2^{1/2}} +\al\)^{2n}=\sum_{j=0}^{2n}  \binom {2n}{ j}    { G^{j}\over 2^{j/2}} \,\,\al^{(2n-j)}.\label{0.3jh}
\end{equation}
This is very suggestive of our Isomorphism Theorem (\ref{1.0}), except that neither the powers $G^{j}$ nor $L_1^{n-j}$ make sense without renormalization. It seems remarkable  that  the subtractions needed for the Wick powers $:G^{j}: ( \nu)$ and those needed for the $L_{n-j}( \nu)$ match up to preserve the simple form of (\ref{1.0}).

Our isomorphism theorem, Theorem \ref{theo-ILT}, for intersection local times of random interlacements has a strong resemblance to the isomorphism theorem, \cite[Theorem 4.2]{MRmem}, for intersection local times $\bar{L}_{n}(\nu)$ of a L\'evy process killed at the end of an independent exponential time $\la$ which states that
\bea
&&
\Big\{\sum_{j=0}^{n}{ n\choose j}\({ :G^{2j}:\over 2^{j}}\times \bar{L}_{n-j}\)(\nu),\,\nu\in D, f(X_{\la})P^{\rho} \times P_{G}\Big\}\label{mem1.0}\\
&&\hspace{1 in}\stackrel{law}{=}
\Big\{{ :G^{2n}:(\nu)\over 2^{n}},\,\nu\in D,   G(\rho)G(f\cdot dx)\, P_{G} \Big\}.\nn
\eea
This holds for any measure $\rho\in  \mathcal{G}_{K}^{1}$ and `nice' function $f$, where
$f(X_{\la})P^{\rho}(F)=\int P^{x}(F\,f(X_{\la}))\,d\rho(x)$. 
Once again, the $n=1$ case of this theorem, essentially due to Dynkin,    suggests an equivalence in law  between
\begin{equation}
G^{2}/2 +\bar L_{1}\hspace{ .2in}\mbox{and}\hspace{ .2in}G^{2}/2\label{0.5jh}
\end{equation}
under the appropriate measures. As before, taking the $n$'th power of both sides formally leads to (\ref{mem1.0}), modulo the renormalizations needed for the terms to make sense.  We do not see how to exploit the formal resemblance between (\ref{mem1.0}) and our 
 isomorphism theorem, Theorem \ref{theo-ILT}. directly, so we must proceed from scratch.
 
 In Section \ref{sec-prelim} we recall various properties of quasi-processes, random interlacements and Wick powers. The proof of our main Theorems require detailed estimates and involved combinatorics. Before proceeding with this, in Section \ref{sec-toy} we show how to obtain the analogue of our isomorphism theorem, Theorem \ref{theo-ILT}, in
a toy model: interlacements of a random walk in $Z^{d}$. Since local times exist in this model, there is no problem with the existence of  renormalized intersection local times, and the combinatorics 
needed for the proof of Theorem \ref{theo-ILT} is greatly simplified.
  Section \ref{sec-rilt} defines the renormalized intersection local times for random interlacements and proves Theorem \ref{theo-multiriltintro-m}. We are able to use many of the results and techniques from \cite{LMR2}, so we only point out the necessary changes. 
  
 The proof of our isomorphism theorem is quite complicated.  
 In Section \ref{sec-gen} we give the proof of our isomorphism theorem, subject to   Lemma \ref{lem-comb} which is proven in Section \ref{sec-comb}. In these sections we are able to use many estimates and combinatorial arguments from  \cite{MRmem}. Again, we concentrate on the differences. Our main effort is in Section \ref{sec-comb} which requires very different combinatorics. 
 
 We mention that \cite{MRmem} is proven under assumptions that are somewhat different than those of the present paper, which is modeled on \cite{LMR2}. However, the reader will have no trouble using the estimates of \cite[Section 8]{LMR2} in place of those used in \cite{MRmem}, so we will freely use the results of  \cite{MRmem}.
 
  The   heuristic formula (\ref{3})   
 involves both self-intersections of paths in the Poisson process $\mathcal{I}_{\al}$ and  intersections between different paths in $\mathcal{I}_{\al}$.  In Section \ref{sec-dec}   we show how to make this explicit.

 Acknowledgements: We are particularly grateful to Kevin O'Bryant for providing the proof of the critical Lemma \ref{lem-comba}, and to Pat Fitzsimmons for discussions about quasi-processes.

\section{Preliminaries on interlacements and Wick powers}\label{sec-prelim}

If $L^{\nu}_{t}, t\geq 0$ denotes the CAF on $\Omega$ with Revuz measure $\nu$ then there is an extension to $W$, which we also denote by $L^{\nu}_{t}, t\in R^{1}$  with the property that
  \begin{equation}
\mu_{m}\(\int_{R^{1}}f\(Y_{t}\)\,dL^{\nu}_{t}\)=\nu (f),\label{q.3}
 \end{equation}
for all measurable $f$, see \cite[XIX, (26.5)]{DM4}. Note in particular that for any bounded  compactly supported measurable  function $g$
\begin{equation}
L^{g\,dm}_{t}=\int_{-\ff}^{ t}g\left( Y_{s}  \right)\,ds.\label{simple.1}
\end{equation}
 
 \bl \label{lem-qpm} For any $\nu_{1},\cdots,\nu_{k}$, with support in some compact $K\subset R^{d}$
\be
\mu_{m}\(\prod_{j=1}^{k}L_{\ff}^{\nu_{j}}\) = \sum_{\pi\in \mbox{\footnotesize Perm$([1,k])$}} \int \prod_{j=1}^{k-1} u(y_{j},y_{j+1}) \prod_{j=1}^{k}\,d\nu_{\pi (j)}(y), \label{q.4}
\ee
where $\mbox{\footnotesize Perm$([1,k])$}$ denotes the set of permutations of $[1,k]$.
\el

{\bf  Proof:  } Let  $T_{K}$ denote the first hitting time of $K$. Then since the $L_{t}^{\nu_{j}}$ do not begin to grow until time $T_{K}$
\begin{eqnarray}
&&\mu_{m}\(\int_{\{-\ff<t_{1}\leq \cdots \leq t_{k-1}\leq t_{k} < \ff\}}\prod_{j=1}^{k}  \,dL_{t_{j}}^{\nu_{j}}\)
\label{q.5}\\
&&=\mu_{m}\(\int_{\{0\leq t_{1}\leq \cdots \leq t_{k-1}\leq t_{k} < \ff\}}\prod_{j=1}^{k}  \,dL_{T_{K}+t_{j}}^{\nu_{j}}\).\nn   \nonumber
\end{eqnarray}
Hence by the second property of $\mu_{m}$
this equals
\begin{equation}
\mu_{m}\(\int_{0}^{\ff} h\(Y_{T_{K}+t_{1}}\)\,dL_{T_{K}+t_{1}}^{\nu_{1}}\),\label{q.6}
\end{equation} 
where
\bea
h(x)&=&E^{x}\(\int_{\{0\leq t_{2}\leq \cdots \leq t_{k-1}\leq t_{k} < \ff\}}\prod_{j=2}^{k}  \,dL_{t_{j}}^{\nu_{j}}\)\label{q.7}\\
&=&\int u(x,y_{2}) \prod_{j=2}^{k-1} u(y_{j},y_{j+1}) \prod_{j=2}^{k}\,d\nu_{j}(y).\nn
\eea
Hence, using  once again  the fact that the $L_{t}^{\nu_{j}}$ do not begin to grow until time $T_{K}$ and then (\ref{q.3}) and
\begin{eqnarray}
&&\mu_{m}\(\int_{\{-\ff<t_{1}\leq \cdots \leq t_{k-1}\leq t_{k} < \ff\}}\prod_{j=1}^{k}  \,dL_{t_{j}}^{\nu_{j}}\)
\label{q.5a}\\
&&=\mu_{m}\(\int_{R^{1}}h\(Y_{t_{1}}\)\,dL^{\nu_{1}}_{t_{1}}\)\nn\\
&&= \int \prod_{j=1}^{k-1} u(y_{j},y_{j+1}) \prod_{j=1}^{k}\,d\nu_{j}(y),\nn   \nonumber
\end{eqnarray}
and (\ref{q.4}) follows, since, up to sets of Lebesgue measure zero \[R^{ k}=\sum_{\pi\in \mbox{\footnotesize Perm$([1,k])$}}\{-\ff<t_{\pi( 1)}\leq \cdots \leq t_{\pi(k- 1)}\leq t_{\pi( k)} < \ff\}.\]
\qed

In particular (\ref{q.4}) shows that
\begin{equation}
\mu_{m}\(\(L_{\ff}^{\nu}\)^{k}\)=k!\int  \prod_{j=1}^{k-1} u(y_{j},y_{j+1})\prod_{j=1}^{k}\nu(dy_{j}).\label{q.8}
\end{equation}

For bounded compactly supported functions $g_{i}$, $ i=1,\ldots,k$, on $R^{d}$, (\ref{q.4}) and (\ref{simple.1})  show that
\bea
\lefteqn{
\mu_{m}\(   \prod_{i=1}^{k}\int_{-\ff}^{\ff}g_{i}(Y_{t})\,dt\) \label{ls.1}}\\
&& = \sum_{\pi\in \mbox{\footnotesize Perm$([1,k])$}}\int u(y_{\pi(1)},y_{\pi(2)})\cdots   u(y_{\pi(k-1)},y_{\pi(k)})  \prod_{i=1}^{k} g_{i}(y_{i})\,dm(y_{i}).\nn
\eea
Simply take $\nu_{i}(dx)=g_{i}(x)\,dm(x)$, $i=1,\ldots,k$.

The potential density $u(x)$ for a symmetric L\'evy process is always weakly positive definite, so there exists a mean zero  Gaussian field $G(\nu)$ for all $\nu\in \GG^{1}$ with covariance 
\begin{equation}
E(G(\nu)G(\nu'))=\int \int u(x-y) \,d\nu(x)\,d\nu'(y).\label{gauss1}
\end{equation}
Since we are interested in the case when $u(0)=\ff$, there is no mean zero  Gaussian process $G_{x}$ with covariance $u(x-y)$. In order to define a substitute for powers of $G$ we proceed as follows.   
 
    Let $f(y)$ be a positive smooth function supported in the unit ball of $R^{d}$ with $\int f(x)\,dx=1$. Set  $f_{\ep}(y)=\ep^{-d}f(y/\ep)$,  and $f_{ \ep,x}(y)=f_{\ep}(y-x)$. 
It follows from our assumptions on $X$ that 
\begin{equation}
u_{\ep,\ep'}( x,x'):=\int u( y,y')\,f_{ \ep,x}(y)\,f_{ \ep',x'}(y')\,dy\,dy'<\ff\label{pd.1}
\end{equation}
 for   $\ep,\ep'>0$. Hence $f_{ \ep,x}(y)\,dy\in \GG^{1}$   and    if we set $G_{x,\ep}=G(f_{ \ep,x}(y)\,dy )$ it follows that
 \begin{equation}
E\left( G_{x,\ep}\,\,G_{x',\ep'} \right)=u_{\ep ,\ep'}( x,x').\label{pd.2}
 \end{equation}
 For $\nu\in  \mathcal{G}^{2}$ we let $:G^{2}:(\nu)$ denote the Wick square corresponding to $\nu$,
a particular   second order Gaussian chaos defined as   
    \begin{equation}
  :G^{2}: (\nu)=\lim_{\ep\rar0} \int \( G^{2}_{x,\ep}-E\(G^{2}_{x,\ep}\)   \) \,d\nu(x). \label{a.3a}
    \end{equation}
(See \cite{MR96}  for details, as well as \cite[Lemma 3.3]{MRmem} for this and the analogue for the higher order Wick powers $:G^{n}: (\nu):$). 
 
 For the proof of the next Theorem only it will be convenient to use another approximation to $G$. It follows from our assumptions that $X$ has symmetric and positive definite transition densities $p_{t}( x,y)$ with potential densities
 \begin{equation}
 u( x,y)=\int_{0}^{ \ff}p_{t}( x,y)\,dt,\label{dep0}
 \end{equation}
  and that
 \begin{equation}
u_{\ep}( x,y):=\int_{\ep}^{ \ff}p_{t}( x,y)\,dt<\ff\label{dep1}
 \end{equation}
 for each $x,y$ and $\ep>0$. Hence $p_{\ep}( x,y)\,dy\in \GG^{1}$   and    if we set $\bar G_{x,\ep}=G(p_{\ep}( x,y)\,dy )$ it follows that
 \begin{equation}
E\left(\bar G_{x,\ep}\,\,\bar G_{x',\ep'} \right)=u_{\ep +\ep'}( x,x').\label{dep.2}
 \end{equation}
If $\nu\in  \mathcal{G}^{2}$ then 
    \begin{equation}
  :G^{2}: (\nu)=\lim_{\ep\rar0} \int \(\bar G^{2}_{x,\ep}-E\(\bar G^{2}_{x,\ep}\)   \) \,d\nu(x). \label{dep.3a}
    \end{equation}
    with convergence in all $L^{ p}$. $\bar G_{x,\ep}$ is convenient since $u_{\ep}( x,y)\uparrow u( x,y) $. However, in the sequel, because of  interlacements it will be important to work with compactly supported $f_{\ep}$ as in the previous paragraph.
 
Set
\begin{equation}
  L_{1}(\nu):=\sum_{\om \in \mathcal{I}_{\al}}L^\nu_\ff(\om).\label{q2.1}
\end{equation}
The next result which  will follow from the master formula for Poisson processes, is a generalization of the Isomorphism Theorems of Sznitman, \cite{Sz2, Sz3,Sz1}. We let  $|\nu|$ denote the mass of $\nu$.

  \bt\label{theo-CAF}
For any  $\al >0$, compact $K\subset S$ and countable  $D\subseteq  \mathcal{G}_{K}^{2}$,
\bea
&&
\Big\{\textstyle{ 1\over 2}:G^{2}:(\nu)+ L_{1}(\nu),\,\nu\in D, \mathcal{P}_{\al^{2}}\times P_{G}\Big\}\label{0.0}\\
&&\stackrel{law}{=}
\Big\{\textstyle{ 1\over 2}:G^{2}:(\nu)+\sqrt{2}\al G(\nu)+\al^{2}  |\nu|,\,\nu\in D,P_{G} \Big\}.\nn
\eea
\et

{\bf  Proof of Theorem \ref{theo-CAF}: }Because everything is additive in $\nu$ we can write (\ref{0.0}) as 
\bea
&&
\mathcal{P}_{\al^{2}}\times P_{G}\(\exp \(\de L_{1}(\nu)+{\de\over 2}:G^{2}:(\nu)\)\)\label{q4.1}\\
&&= P_{G}\(\exp \(\textstyle{ \de\over 2}:G^{2}:(\nu)+\de\sqrt{2}\al G(\nu)+\de\al^{2}  |\nu|\)\)\nn
\eea
for $\de$ small. Equivalently, we show that
\be
\mathcal{P}_{\al^{2}}\(e^{\de L_{1}(\nu)}\)= { P_{G}\(\exp \(\textstyle{ \de\over 2}:G^{2}:(\nu)+\de\sqrt{2}\al G(\nu)+\de\al^{2}  |\nu|\)\)\over P_{G}\(\exp \({\de\over 2}:G^{2}:(\nu)\)\)}.\label{q4.2}
\ee

We first note that using (\ref{a.3a}), the Gaussian moment formula and the monotone convergence theorem  we have 
\begin{equation}
P_{G}\(\exp \({\de\over 2}:G^{2}:(\nu)\)\)=\lim_{\ep\to 0}P_{G}\(\exp \({ \de\over 2}\int \(\bar  G^{2}_{x,\ep}-E\(\bar G^{2}_{x,\ep}\)   \) \,d\nu(x)\)\)\label{q4.4}
\end{equation}
and
\bea
&&
P_{G}\(\exp \(\textstyle{ \de\over 2}:G^{2}:(\nu)+\de\sqrt{2}\al G(\nu)+\de\al^{2}  |\nu|\)\)\label{q4.5}\\
&&=\lim_{\ep\to 0}P_{G}\(\exp \({ \de\over 2}\int \((\bar G_{x,\ep}+\sqrt{2}\al)^{2} -E\(\bar G^{2}_{x,\ep}\)   \) \,d\nu(x)\)\).\nn
\eea
Therefore, 
\begin{eqnarray}
&&  { P_{G}\(\exp \(\textstyle{ \de\over 2}:G^{2}:(\nu)+\de\sqrt{2}\al G(\nu)+\de\al^{2}  |\nu|\)\)\over P_{G}\(\exp \({\de\over 2}:G^{2}:(\nu)\)\)}
\label{q4.6}\\
&&=\lim_{\ep\to 0}   {P_{G}\(\exp \({ \de\over 2}\int \((\bar G_{x,\ep}+\sqrt{2}\al)^{2} -E\(\bar G^{2}_{x,\ep}\)   \) \,d\nu(x)\)\) \over P_{G}\(\exp \({ \de\over 2}\int \(\bar  G^{2}_{x,\ep}-E\(\bar G^{2}_{x,\ep}\)   \) \,d\nu(x)\)\)}\nonumber\\
&&=\lim_{\ep\to 0}   {P_{G}\(\exp \({ \de\over 2}\int  (\bar G_{x,\ep}+\sqrt{2}\al)^{2}      \,d\nu(x)\)\) \over P_{ G}\(\exp \({ \de\over 2}\int  \bar  G^{2}_{x,\ep}     \,d\nu(x)\)\)}.\nonumber
\eea
A simple Gaussian computation, see \cite[Lemma 5.2.1]{book}, shows that  this is 
\bea
&&=\lim_{\ep\to 0}  \exp \(\al^{2} \(\sum_{n=1}^{\ff}\de^{n }\int \prod_{j=1}^{n-1} u_{2\ep}(x_{j},x_{j+1})\prod_{j=1}^{n}\nu(dx_{j})\)\) \nonumber\\
&&=\exp \(\al^{2} \(\sum_{n=1}^{\ff}\de^{n }\int \prod_{j=1}^{n-1} u(x_{j},x_{j+1})\prod_{j=1}^{n}\nu(dx_{j})\)\).\nn
\end{eqnarray}
by the monotone convergence theorem.

On the other hand, by the master formula for Poisson processes, \cite{K},
\begin{equation}
\mathcal{P}_{\al^{2}}\(e^{\de  L_{1}(\nu)}\)=\exp \(\al^{2}\mu_{m}\(e^{\de L^\nu_{\ff}}-1\)\),\label{4.mast}
\end{equation}
and it follows from (\ref{q.8}) that 
\be
\mu_{m}\(e^{\de L^\nu_{\ff}}-1\)=\sum_{n=1}^{\ff}\de^{n } \,\,{\mu_{m}\((L^\nu_{\ff})^{n}\)\over n!} =\sum_{n=1}^{\ff}\de^{n }\int \prod_{j=1}^{n-1} u (x_{j},x_{j+1})\prod_{j=1}^{n}\nu(dx_{j}).
\label{4.fm}
\ee
This completes the proof of (\ref{q4.2}) and hence of (\ref{0.0}).\qed

   The  next lemma, which describes the moment structure of random interlacements, 
  follows from (\ref{ls.1}) and  the master formula for Poisson processes. It
will be used in the next section. For any bounded compactly supported function $g$ on $R^{d}$
we set
\begin{equation}
L_{1}(g):=L_{1}(g(x)\,dx).\label{deff}
\end{equation}
Then by (\ref{q2.1}) and (\ref{simple.1})
\begin{equation}
  L_{1}(g)=\sum_{\om \in \mathcal{I}_{\al}}\int g(Y_{t}(\om))\,dt.\label{simple.2}
\end{equation}

\begin{lemma}\label{lem-rim} Let   $g_{j}$, $ j=1,\ldots,k$ be bounded compactly supported functions on $R^{d}$. Then  
\bea
&&
\mathcal{P}_{\al}\(   \prod_{i=1}^{k}L_{1} ( g_{i})\)\label{ls.3}\\
&&=\sum_{\stackrel{B_{1}\cup \cdots\cup B_{j}=[1,k]}{j=1,\ldots,k}}\al^{j}\prod_{l=1}^{j} \mu_{m}\(   \prod_{i\in B_{l}} \int_{-\ff}^{\ff}g_{i}(Y_{t})\,dt\)\nn\\
&& \hspace{-.3 in}=\sum_{\stackrel{B_{1}\cup \cdots\cup B_{j}=[1,k]}{j=1,\ldots,k}}\al^{j} \int\(\prod_{l=1}^{j}
\sum_{\pi\in \mbox{\footnotesize Perm$(B_{l})$}} u(y_{\pi(1_{l})},y_{\pi(2_{l})})\cdots   u(y_{\pi((|B_{l}|-1)_{l})},y_{\pi(|B_{l}|_{l})}) \)\nn\\
&&\hspace{3.8in} \prod_{i=1}^{k} g_{i}(y_{i})\,dy_{i}\nn
\eea
where the sum in the second and third line is over all partitions of $[1,k]$ and  $\mbox{\footnotesize Perm$(B_{l})$}$ denotes the set of permutations of $B_{l}=\{1_{l}, 2_{l},\ldots, |B_{l}|_{l}\}$.
\end{lemma}

\section{A toy model: interlacements of random walks in $Z^{d}$}\label{sec-toy}

In this section we take $X$ to be a symmetric continuous time transient random walk in $Z^{d}$. We take $m$
 to be counting measure, and as before we denote the potential as $u(x,y)=u(x-y)$. Here $u$ is finite.  Let
 \begin{equation}
  L_{1}(y):=\sum_{\om \in \mathcal{I}_{\al}}L^y_\ff(\om)\label{rw.1}
\end{equation}
 where $L^y_\ff(\om)$ is the total local time of the path $\om$ at $y\in Z^{d}$. We refer to this as a toy model for intersections because when local times exist, intersection local times are straightforward, and in particular there is no real need for renormalization. Nevertheless, it gives us an opportunity to exhibit some of  the combinatorics involved in proving our main theorem without the need to 
 to deal with approximations and error bounds.
 
 In the following we abbreviate $u=u(0)$. Then if $\{G_{x}, x\in Z^{d}\}$ denotes the Gaussian process with covariance $u(x,y)$, the Wick powers $:G_{x}^{n}: $ are defined as 
 \begin{equation}
:G_{x}^{n}:=\sum_{j=0}^{[n/2]} (-1)^{j}{n \choose 2j} {(2j)! \over j!2^{j}}u^{j}G_{x}^{n-2j}.     \label{rw.2}
 \end{equation}
 Thus $:G_{x}^{n}: $ is an $n$'th degree polynomial in $G_{x}$, and it has generating function
  \begin{equation}
  \sum_{n=0}^{\ff}{s^{n}:G_{x}^{n}: \over n!}  =e^{sG_{x}-s^{2}u/2}.\label{15.6g}
  \end{equation}
  See \cite[(3.8), (3.16)]{MRmem}. We will use the convention that if $f(x)=\sum_{n=0}^{\ff}a_{n}z^{n}$
   is analytic, then $:f(G_{x}):=\sum_{n=0}^{\ff}a_{n}:G_{x}^{n}:$
  
Since $:G_{x}^{2}:=G_{x}^{2}-u, $ Theorem \ref{theo-CAF}  for interlacements of random walks takes the form:
for any  $\al >0$, compact $K\subset Z^{d}$,
\bea
&&
\Big\{\textstyle{ 1\over 2}G_{x}^{2} + L_{1}(x),\,x\in K, \mathcal{P}_{\al^{2}}\times P_{G}\Big\}\label{rw2.0}\\
&&\stackrel{law}{=}
\Big\{\textstyle{ 1\over 2}G_{x}^{2}+\sqrt{2}\al G_{x}+\al^{2},\,x\in K,P_{G} \Big\}.\nn
\eea
It is interesting to compare this with the generalized second Ray-Knight Theorem, \cite{five}.
 
The goal of this section is to prove Theorem \ref{theo-ILT} for random walks: 
\bea
&&
\Big\{\sum_{j=0}^{n}{ n\choose j}\({ :G_{x}^{2j}:\over 2^{j}}  L_{n-j}(x)\),\,x\in K, \mathcal{P}_{\al^{2}}\times P_{G}\Big\}\label{rwilt.0}\\
&&\stackrel{law}{=}
\Big\{\sum_{j=0}^{2n}{ 2n\choose j}\al^{(2n-j)}{ :G_{x}^{j}:\over 2^{j/2}},\,x\in K,P_{G} \Big\}.\nn
\eea
 Here $L_{0}(x)=1$ and  $L_{n}(x)$ for $n>1$ is defined by
    \begin{equation}
  \sum_{n=0}^{\ff}{s^{n} L_{n}(x)\over n!} =e^{{s\over 1+su}L_{1}(x) }.\label{15.2b}
  \end{equation}
  
Note that
\begin{equation}
\(G_{x}/\sqrt{2}+\al\)^{2}= { 1 \over 2}G_{x}^{2}+\sqrt{2}\al\,G_{x}+\al^{2}.\label{rwH1}
\end{equation}

  Let $\MM(\JJ_{K})$ denote the set of functions measurable with respect to
$\mathcal {J}_{K}=: \si\(\(G_{x}/\sqrt{2}+\al\)^{2};\,x\in K\)$. We define the ring
homomorphism \be
\Phi:\MM(\JJ_{K})\mapsto \MM(\JJ_{K}\times \FF)\label{1mmc} \ee 
as the measurable
extension of the mapping $\Phi$ such that $\Phi(1)=1$ and
\be
\Phi\(\prod_{i=1}^n \(G_{x_{i}}/\sqrt{2}+\al\)^{2}\)=\prod_{i=1}^n \(   { G_{x_{i}}^{2}\over 2}+L_1(x_i) \),
\hspace{.2in}n=1,\ldots,\label{rwphi1}
\ee 
where $\FF$ is the $\si$-algebra generated by the random interlacement.  With this
notation  
(\ref{rw2.0}) can be reformulated as follows: Let $(h_1,h_2,\ldots)$ be a sequence of
$\JJ_{K}$ measurable functions. Then for any Borel measurable non-negative
function $F$ on $R^{\ff}$
\be E_{G}\mathcal{P}_{\al^{2}}\(F(\Phi(h_1),\Phi(h_2),\ldots)\)
=E_{G}\(F(h_1,h_2,\ldots) \).\label{c15.3j}
\ee 
Hence to prove   (\ref{rwilt.0}) it suffices to show that if we define
\begin{equation}
J_{n}(x):=:\(G_{x}/\sqrt{2}+\al\)^{2n}:=\sum_{j=0}^{2n}{ 2n\choose j}\al^{(2n-j)}{ :G_{x}^{j}:\over 2^{j/2}}\,, \label{15.20s}
\end{equation}
then  
\begin{equation}
J_{n}(x)\in  \si\(\(G_{x}/\sqrt{2}+\al\)^{2}\),\label{15.20}
\end{equation}  
and
\begin{equation}
\Phi\( J_{n}(x)\)=\sum_{j=0}^{n}{ n\choose j}\({ :G_{x}^{2j}:\over 2^{j}}  L_{n-j}(x)\).\label{15.21}
\end{equation}

 We abbreviate $G=G_{x}, L_{n}=L_{n}(x)$.
The following Lemma is proven below.
\bl\label{lem-lag}
\begin{equation}
  \sum_{n=0}^{\ff}{s^{n}:\(G/\sqrt{2}\)^{2n}:\over n!}=(1+us)^{-1/2}\exp \({s (G/\sqrt{2})^{2}\over 1+us}\),\label{lag.1}
\end{equation}
and
\begin{equation}
  \sum_{n=0}^{\ff}{s^{n}:\(G/\sqrt{2}+\al\)^{2n}:\over n!} =(1+us)^{-1/2}\exp \({s (G/\sqrt{2}+\al)^{2}\over 1+us}\).\label{lag.2}
\end{equation}
\el

(\ref{15.20}) follows from (\ref{lag.2}), and then applying  $\Phi $ to both sides of (\ref{lag.2}) we obtain
        \bea
  \sum_{n=0}^{\ff}{s^{n}\Phi \(:\(G/\sqrt{2}+\al\)^{2n}:\)\over n!} 
   &=&(1+us)^{-1/2}\exp \({s (G^{2}/2+L_{1})\over 1+us}\)\nn\\
   &=&(1+us)^{-1/2}\exp \({s  G^{2}/2 \over 1+us}\)\exp \({s  L_{1} \over 1+us}\)\nn
  \eea
  so that by (\ref{lag.1}) and (\ref{15.2b})
  \begin{eqnarray}
  &&  \sum_{n=0}^{\ff}{s^{n}\Phi \(:\(G/\sqrt{2}+\al\)^{2n}:\)\over n!}
  \label{}\\
  &&=\(   \sum_{m=0}^{\ff}{s^{m}:\(G/\sqrt{2}\)^{2m}:\over m!} \)  \( \sum_{j=0}^{\ff}{s^{j} L_{j}\over j!}\)\nonumber
  \end{eqnarray}
  which easily proves (\ref{15.21}).

{\bf  Proof of Lemma \ref{lem-lag}: }
By \cite[8.957.1]{GR}, the generating function for the Hermite polynomials $H_{n}$ is
    \begin{equation}
  \sum_{n=0}^{\ff}{z^{n}H_{n}(x) \over n!} =e^{2tx-t^{2}}.\label{15.7k}
  \end{equation} 
  Setting $x=G/\sqrt{2u}, z=s\sqrt{u/2}$ and comparing with (\ref{15.6g}) we see that
  \begin{equation}
 :G^{n}:=(u/2)^{n/2} H_{n}\(G/\sqrt{2u}\).\label{15.8}
  \end{equation}
  If we use the notation $\mathcal{L}_{n}^{\ka}(x)$ for the $n$'th order Laguerre polynomial of index $\ka$, it follows from  \cite[8.972.2]{GR} that
  \begin{equation}
  H_{2n}(x)=(-1)^{n}n!2^{2n}\mathcal{L}_{n}^{-1/2}(x^{2}),\label{15.9}
  \end{equation}
  hence using the previous formula
    \begin{equation}
{ :G^{2n}: \over 2^{n}}=  (-u)^{n}n! \mathcal{L}_{n}^{-1/2}(G^{2}/2u).\label{15.10}
  \end{equation}
Therefore
     \bea
  \sum_{n=0}^{\ff}{s^{n}:\(G/\sqrt{2}\)^{2n}:\over n!} &=&  \sum_{n=0}^{\ff}(-us)^{n}\mathcal{L}_{n}^{-1/2}(G^{2}/2u)\label{15.27}\\&=& (1+us)^{-1/2}\exp \({s G^{2}/2\over 1+us}\)\nn
  \eea
  by \cite[8.975.1]{GR}. This gives (\ref{lag.1}).

  By (\ref{rw.2}), or (\ref{15.10})
  \begin{equation}
 :G^{2n}:=P_{n}(G^{2}) \label{15.31}
  \end{equation}
  for some $n$'th degree polynomial $P_{n}$.
In view of (\ref{lag.1}), (\ref{lag.2}) is equivalent to the following: for any $c\in R^{1}$
  \begin{equation}
 :(G+c)^{2n}:=P_{n}((G+c)^{2}) \label{15.32}
  \end{equation}
  for the same polynomial $P_{n}$.

But by (\ref{15.6g})
  \begin{equation}
\(  \sum_{i=0}^{\ff}{s^{i}:G^{i}: \over i!} \)\(  \sum_{j=0}^{\ff}{s^{j}c^{n} \over j!} \) =e^{sG-s^{2}u/2}e^{sc},\label{15.46}
  \end{equation}
  or equivalently
    \begin{equation}
   \sum_{n=0}^{\ff}{s^{n} :(G+c)^{n}:\over n!}   =e^{s(G+c)-s^{2}u/2},\label{15.46a}
  \end{equation}
which establishes (\ref{15.32}).

\section{Renormalized intersection local times}\label{sec-rilt}

Let 
\begin{equation}
L_{1}(x,\ep):=L_{1} ( f_{\ep,x} )=\sum_{\om \in \mathcal{I}_{\al}}\int  f_{\ep,x} (Y_{t}(\om))\,dt. \label{rilt.1q}
\end{equation} 
  $L_{1}(x,\ep)$ can be thought of as the approximate total  local time of the random interlacement at the point $x\in R^{d}$.  When $u(0)=\ff$,  local times do  not exist and we can not take the limit of $L_{1}(x,\ep)$ as $\ep\to 0$.    Nevertheless, it is often the case that renormalized intersection local times exist.    We proceed to define renormalized intersection local times.
  
 We begin with the definition of the chain functions 
  \begin{equation}
  \mbox{ch}_{k}(r) =  \int   u (ry_{1},ry_{2}) \cdots u(ry_{k},ry_{k+1})   \prod_{j=1}^{k+1} f(y_{j}) \,dy_{j},\hspace{.2 in} k\geq 1.\label{1L.1}
  \end{equation} 
  Note that $\mbox{ch}_{k}(r)$ involves $k$ factors of the potential density $u$,   but $k+1$ variables of integration.
    For any $  \si=(k_{1}, k_{2},\ldots )$ let   
      \begin{equation}
 |\si|=\sum_{i=1}^{\ff}  ik_{i}\qquad\mbox{and}\qquad  |\si|_{+}=\sum_{i=1}^{\ff}   (i+1)k_{i}.\label{rilt.2a}
  \end{equation}
  We   define  recursively 
  \begin{eqnarray}
&&
L_{n}(x,r)  = L_{1}^{n} (x,r)  -\sum_{\{ \si\,|\, 1\leq | \si| < | \si|_{+} \leq n\} } J_{n}( \si,r) \label{rilt.1},
  \end{eqnarray}
  where
\begin{equation}
J_{n}( \si,r)={n! \over \prod_{i=1}^{\ff}k_{i}! (n-| \si|_{+})!}\prod_{i=1}^{\ff} \(\mbox{ch}_{i}( r)\)^{k_{i} } 
    L_{n-|\si|} (x,r) .\label{rilt.2b}
\end{equation} 
(Note that $n-|\si|_{+} \ge 0$.)

To help in understanding (\ref{rilt.1}) we note that 
\begin{equation}
L_{2}(x,r)  = L_{1}^{2} (x,r)  -2 \,  \mbox{ch}_{1}(r)L_{1}(x,r)\label{rilt.2ex}
\end{equation}  
and 
\bea
L_{3}(x,r) & =& L_{1}^{3} (x,r)  -6\,  \mbox{ch}_{1}(r)L_{2}(x,r)-6\mbox{ch}_{2}(r)L_{1}(x,r)\label{rilt.3ex}\\
& =&  L_{1}^{3} (x,r)  -6\,   \mbox{ch}_{1}(r)L_{1}^{2}(x,r)+(12\,  \mbox{ch}^{2}_{1}(r)-6\, \mbox{ch}_{2}(r))L_{1}(x,r).\nn
\eea 

It is   interesting to note that one can define $L_{n}(x,r)$ directly, because   
  $L_{n}(x,r)=B_{n}\(L_{1}(x,r)\)$ where the polynomials $B_{n}(u)$ satisfy
  \begin{equation}
  \sum_{n=0}^{\ff}{\(   \sum_{j=0}^{\ff}\mbox{ch}_{j}( r)\,t^{j} \)^{n} \over n!}t^{n}B_{n}(u)=e^{tu},\label{ljo}
  \end{equation}
  where we set $ \mbox{ch}_{0}(r)=1,\,B_{0}(u)=1$. We remark that for our toy model of Section \ref{sec-toy}, $\mbox{ch}_{j}( r)=u^{j}$ and  (\ref{ljo}) takes the form
  \begin{equation}
  \sum_{n=0}^{\ff}{\(  {t\over 1-tu} \)^{n} \over n!}L_{n} =e^{tL_{1} }.\label{15.1}
  \end{equation}
Setting $t=s/(1+su)$ we obtain (\ref{15.2b}).

  The next theorem  gives the joint moments of   the  $L_{n}(\nu)$.  

 \begin{theorem}   \label{theo-1.3}   Let $\mathcal{I}_{\al}$ be the random interlacement   associated with
 a L\'evy process  with potential density $u$ as in Theorem  \ref{theo-multiriltintro-m}. 
Let $n=n_{1}+\cdots+n_{k}$,        and $ \nu_{i}\in \mathcal{G}_{K}^{2n_{i}}$. Then 
 \bea 
 &&
\mathcal{P}_{\al}\(  \prod_{i=1}^{k} L_{n_{i}}(\nu_{i} ) \)\label{rilt.15intro}\\
&&
= \prod_{i=1}^{k}n_{i}! \sum_{B_{1}\cup \cdots\cup B_{j}=[1,n]}\al^{j} \int  \sum_{\pi\in \mathcal{M}_{a}} \prod_{l=1}^{j}  \prod_{i=1}^{|B_{l}|-1}u (x_{\pi(i_{l})},x_{\pi((i+1)_{l})}) \prod_{m=1}^{k}\,d\nu_{m}  (x_{m} )\nn
\eea 
 where the first sum is over all partitions of $[1,n]$, we denote the elements of $B_{l}$ by $\{1_{l}, 2_{l},\ldots, |B_{l}|_{l}\}$ and $\mathcal{M}_{a}$ is the set of maps $\pi:[1,n]\mapsto [1,k] $  with
 $|\pi^{-1}(m)|=n_{m}$ for each $m$  and  such that    for each $l$,  if $\pi(i_{l})=m$  then $\pi((i+1)_{l})\neq m$. (The subscript `a'  in $\mathcal{M}_{a}$ stands for alternating).
 \et

\medskip	\noindent {\bf  Proof  of Theorems   \ref{theo-multiriltintro-m} and  \ref{theo-1.3}  }  The  proof is very similar to the proof of  \cite[Theorem 1.3]{LMR2}.  We begin by showing that for the approximate identities $f_{r,x}$  and bounded functions   $f_{i}=f_{r_{i}, x_{i}},i= 1,\ldots,m$,     
  \bea
&&
\mathcal{P}_{\al}\(L_{n }(x,r)   \prod_{i=1}^{m}L_{1 }(f_{i}) \)\label{20.15s}\\
&&\hspace{-.3 in} =\sum_{\stackrel{B_{1}\cup \cdots\cup B_{j}=[1,m+n]}{j=1,\ldots,m+n}}\al^{j} \nn\\
&&\int\(\prod_{l=1}^{j}
\sum_{\pi_{l}\in \mbox{\footnotesize Perm$_{m,n}(B_{l})$}} u(y_{\pi_{l_{l}}(1)},y_{\pi_{l}(2_{l})})\cdots   u(y_{\pi_{l}((|B_{l}|-1)_{l})},y_{\pi_{l}(|B_{l}|_{l})}) \) \nn\\
&&  \prod_{i=1}^{m} f_{i}  (y_{i} )\,dy_{i} \prod_{i=m+1}^{m+n}f_{r,x}  (y_{i} )\,dy_{i} + \int E_{r}(x,  {\bf z })\prod_{i=1}^{m} f_{i}  (z_{i} )\,dz_{i} \prod_{i=m+1}^{m+n}f_{r,x}  (z_{i} )\,dz_{i} \nn,
\eea
  where 
  $\mbox{\footnotesize Perm$_{m,n}(B_{l})$}$ is the subset of permutations $\pi_{l}$ of $B_{l}=\{1_{l}, 2_{l},\ldots, |B_{l}|_{l}\}$  with the property   that for all $i_{l},(i+1)_{l}\in B_{l}$, if $\pi_{l}(i_{l})\in [m+1,m+n]$, then $\pi_{l}((i+1)_{l})\in [1,m]$.    That is, under the permutation $\pi_{l}$, no two elements of $[m+1,m+n]$ are adjacent.

   The last term in (\ref{20.15s}) is an error term.  It is actually the sum of many terms, some of which may depend on some of the $z_{1},\ldots,z_{n+m}  $.  We use ${\bf z }$ to designate $z_{1},\ldots,z_{n+m}  $.   Since   the  $f$'s are  probability 
density functions we write  last term in (\ref{20.15s})  as an expectation,
 \be 
\EE_{{\bf f }}(E_{r}(x,  {\bf z })):=\int E_{r}(x,  {\bf z })\prod_{i=1}^{m} f_{i}  (z_{i} )\,dz_{i} \prod_{i=m+1}^{m+n}f_{r,x}  (z_{i} )\,dz_{i}.\label{2.13a}
\ee
 We show later that   for   $\nu\in \mathcal{G}_{K}^{2n}$,
  \begin{equation}
 \lim_{r\rar 0} \sup_{\forall |z_{i}|\le 1}\int E_{r}(x, {\bf z }) \,d\nu(x)=0,\label{err0}
  \end{equation}
  which implies that 
  \begin{equation}
    \lim_{r\rar 0}\int  \EE_{{\bf f }}(E_{r}(x,  {\bf z }))\,d\nu(x)=0.
   \end{equation}
(Note that since $f$ is supported on the unit ball in $R^{d}$ we can take $|z_{i}|$ bounded uniformly for all indices $i$.) We deal with all the additional  error terms that are introduced similarly.

\medskip	 Assume that   (\ref{20.15s}) is proved for $L_{n'}(x,r)$, $n'<n$. For any $  \si=(k_{1}, k_{2},\ldots )$ let $\mbox{\footnotesize Perm$_{m+n}(\si)$}$  denote the set of permutations $\bar\pi$  of $[1,m+n]$ of the form $\bar\pi=(\pi_{1},\ldots,\pi_{j})$,  (where  $\pi_{l}$ is a permutation of $B_{l}=\{1_{l}, 2_{l},\ldots, |B_{l}|_{l}\}$ with $B_{1}\cup \cdots\cup B_{j}=[1,m+n]$) 
that contain    $k_{i}$
chains of order $i=1,2,\ldots$ in $[m+1,m+n]$.   (A chain of order $i\geq 1$ is  a  sequence  $ \pi_{l}(j_{l}),\pi_{l}((j+1)_{l}),\ldots,  \pi_{l}((j+i)_{l})$ in $[m+1,m+n]$ for some $l$ which is maximal in the sense that $ \pi_{l}((j-1)_{l}),$ and $ \pi_{l}((j+i+1)_{l})$ are not in $[m+1,m+n]$, possibly because $j-1=0$ or $j+i=|B_{l}|$.)   

Let $\mbox{\footnotesize Perm$_{m,n}$}$  denote the set of permutations $\pi'$  of $[1,m+n]$ of the form $\pi'=(\pi_{1},\ldots,\pi_{j})$,  where each   $\pi_{l}$ is a permutation in \mbox{\footnotesize Perm$_{m,n}(B_{l})$} with $B_{1}\cup \cdots\cup B_{j}=[1,m+n]$.
As in the proof of  \cite[Theorem 1.3]{LMR2}, we see that  the term for any  $\bar\pi\in \mbox{\footnotesize Perm$_{m+n}(\si)$}$ in the evaluation of 
\be
\mathcal{P}_{\al}\(L_{1}^{n}(x,r)   \prod_{i=1}^{m}L_{1 }(f_{i})\),\label{2.20}
 \ee 
is the same as  the term in (\ref{20.15s}) for a particular permutation $\pi'\in \mbox{\footnotesize Perm$_{m,n-|\si|}$}$     in the evaluation of 
\be
\mathcal{P}_{\al}\( \prod_{i=1 }^{\ff} \(\mbox{ch}_{i}( r)\)^{k_{i}} 
L_{n-|\si|}(x,r) \prod_{i=1}^{m}L_{1 }(f_{i})\),  \label{2.21}
   \ee 
   up to error terms  $H_{r}(x,{\bf z })$.   (We note that $\pi'$ will be associated with a partition 
   $B'_{1}\cup \cdots\cup B'_{j}=[1,m+n-|\si|]$, where $B'_{l}\subseteq B_{l}$ for each $l$).
And as in that proof we can do the combinatorics to show that   up to the error terms, the contribution to (\ref{2.20}) from $\mbox{\footnotesize Perm$_{m+n}(\si)$}$ is equal to 
\bea
&&
{n! \over \prod_{i=1}^{\ff}k_{i}! (n-|\si|_{+})!} \mathcal{P}_{\al}\( \prod_{i=1 }^{\ff} \(\mbox{ch}_{i}( r)\)^{k_{i}} 
L_{n-|\si|}(x,r) \prod_{i=1}^{m}L_{1 }(f_{i})\)  \nn\\
&&\qquad=\mathcal{P}_{\al}\( J_{n}(\si,r)\prod_{i=1}^{m}L_{1 }(f_{i})\).\label{2.21s}
   \eea 
   If we let $\mbox{\footnotesize Perm$_{m+n}$}$  denote the set of permutations $\wh\pi$  of $[1,m+n]$ of the form $\wh \pi=(\pi_{1},\ldots,\pi_{j})$,  where each   $\pi_{l}$ is a permutation     of $B_{l}$ with $B_{1}\cup \cdots\cup B_{j}=[1,m+n]$, then
 considering (\ref{rilt.1}) and the fact that 
$\mbox{\footnotesize Perm$_{m+n}$}-\mbox{\footnotesize Perm$_{m,n}$}=\cup_{|\si|\geq 1}\mbox{\footnotesize Perm$_{m+n}(\si)$}$, we see that   the induction step in the proof of (\ref{20.15s}) is proved.

\medskip 	We iterate the steps  used in the proof of   (\ref{20.15s}), and use the fact that each of the $L_{n_{i}}(x_{i}, r  )$ are sums of multiples of $L(x_{i}, r  )$ to obtain   
 \bea
&&
\mathcal{P}_{\al}\(  \prod_{i=1}^{k} L_{n_{i}}(x_{i}, r_{i}  ) \) =\sum_{\stackrel{B_{1}\cup \cdots\cup B_{j}=[1,m+n]}{j=1,\ldots,m+n}}\al^{j}\label{20.15aa}\\
&& \int\(\prod_{l=1}^{j}
\sum_{\pi_{l}\in \mbox{\footnotesize Perm$_{n_{1},\ldots,n_{k}}(B_{l})$}} u(y_{\pi_{l_{l}}(1)},y_{\pi_{l}(2_{l})})\cdots   u(y_{\pi_{l}((|B_{l}|-1)_{l})},y_{\pi_{l}(|B_{l}|_{l})}) \)\nn\\
&&\hspace{1 in}\prod_{j=1}^{n} f _{r_{g(j)},x_{g(j)}}  (y_{j} )\,dy_{j}+\EE_{{\bf f }}(E_{r_{1},\ldots, r_{k}}(x_{1},\ldots, x_{k},{\bf z })), \nn
\eea
where   $\mbox{\footnotesize Perm$_{n_{1},\ldots,n_{k}}(B_{l})$}$ is the set of permutations   $\pi_{l}$  of $B_{l}=\{1_{l}, 2_{l},\ldots, |B_{l}|_{l}\}$  with the property   that for all $m$, when $\pi_{l}(m_{l})\in \big[1+\sum_{p=1}^{i-1}n_{p},\sum_{p=1}^{i }n_{p}\big]:=C_{i}$ then $\pi_{l}((m+1)_{l})\notin C_{i}$,  for all $i\in [1,k]$,  and     $g(m_{l})=i$ when  $m_{l}\in C_{i}$. (In the last term in (\ref{20.15aa}) we use the notation introduced in (\ref{2.13a}).)

The error terms are very similar to those  described in the proof of  \cite[Theorem 1.3]{LMR2}, the only difference being  that there we worked with the loop measure whereas here we have the quasi-process measure $\mu_{m}$.  

Set  
\be L_{n,r}(\nu)=  \int   L_{n}(x,r) \,d\nu(x).\label{2.47}
\ee
Following the proof of  \cite[Theorem 1.3]{LMR2} we can then show that
 of Theorem \ref{theo-multiriltintro-m}, and then  of Theorem \ref{theo-1.3}.\qed

\section{Proof of the Isomorphism Theorem}\label{sec-gen}

We  define the     cycle functions
    \begin{equation}
 \mbox{cy}_{k}(\ep)  =    \int   u (\ep y_{1},\ep y_{2}) \cdots u (\ep y_{k-1},\ep y_{k})u (\ep y_{k},\ep y_{1})  \prod_{j=1}^{k} f(y_{j}) \,dy_{j}.\label{4.2}
  \end{equation} 
  For any $\si=(k_{1}, k_{2},\ldots ; m_{2}, m_{3},\ldots)$ we set 
  \begin{equation}
 |\si|=\sum_{i=1}^{\ff}  ik_{i}+\sum_{j=2}^{\ff}  jm_{j},\hspace{.2 in}  |\si|_{+}=\sum_{i=1}^{\ff}   (i+1)k_{i}+\sum_{j=2}^{\ff}  jm_{j}.\label{4.3}
  \end{equation}
  For example, $\si=(2;0,1)$ means that $k_{1}=2, m_{3}=1$ and all other $k_{i}, m_{j}=0$, and 
 $ |\si|=5,  |\si|_{+1}=7.$

Set $\wt H_{0}(x,\ep)  =1 $ and $ \wt H_{1}(x,\ep)  = H_{1} (x,\ep)  ={ :G^{2}:(x,\ep) \over 2} +2^{1/2}\al\, G(x,\ep)  +\al^{2}$.
We then  define inductively
  \be
\wt H_{n}(x,\ep)   =H_{1}^{n} (x,\ep)   -\sum_{\{\si\,|\, 1\leq |\si| \leq |\si|_{+} \leq n\} } I_{n, \ep }( \si)\wt H_{n-|\si| } (x,\ep) , \label{4.4}
  \ee
  where
\begin{equation}
I_{n, \ep }( \si)={n! \over (n-|\si|_{+})!\prod_{i=1}^{\ff}k_{i}!\prod_{j=2}^{\ff}m_{j}!}\prod_{i=1}^{\ff} \(\mbox{ch}_{i}(\ep )\)^{k_{i}} 
  \prod_{j=2}^{\ff} \({\mbox{cy}_{j}(\ep ) \over 2j}\)^{m_{j}}.\label{4.5}
\end{equation} 

For any $\al>0$ let
 \begin{equation}
H_{n}(\nu):= \sum_{j=0}^{2n}{ 2n\choose j}\al^{(2n-j)}{ :G^{j}:(\nu)\over 2^{j/2}}.\label{4.0}
 \end{equation}
The following is the key technical result of this paper.

\bl\label{lem-comb}For any $\nu\in\mathcal{G}_{K}^{2n'}$ and $\al\geq 0$
\begin{equation}
\lim_{\ep\to 0}\int \wt H_{n}(x,\ep) \,d\nu(x)=H_{n}(\nu),\label{com.1}
\end{equation}
for all $1\leq n\leq n'$.
\el
 
When $\al=0$ this is \cite[Lemma 4.3]{MRmem}.

 It is easy to check that
 \bea
\wt H_{2}(x,\ep)   &=&H_{1}^{2} (x,\ep) -I_{2,\ep }(1;0)H_{1}  (x,\ep) -I_{2,\ep }(0;1)\label{4.6}\\
 &=&H_{1}^{2} (x,\ep) -2 \mbox{ch}_{1}( \ep )H_{1}  (x,\ep) - {\mbox{cy}_{2}(\ep)  \over 2},\nn 
 \eea 
compare (\ref{rilt.2ex}).

Lemma \ref{lem-comb} will be proven in the next section. Now we return to $L^{n}(x,\ep) $ in order to establish our Isomorhism Theorem.

It follows from (\ref{rilt.1})-(\ref{rilt.2b}) that
  \be 
L_{1}^{n}(x,\ep)  =\sum_{\{\bar\si\,|\, 0\leq |\bar\si| \leq |\bar\si|_{+} \leq n\} } {n! \over \prod_{i=1}^{\ff}k_{i}! }\prod_{i=1}^{\ff} \(\mbox{ch}_{i}( \ep)\)^{k_{i}} 
   {  L_{n-|\bar\si|} (x,\ep)\over (n-|\bar\si|_{+})!}.\label{11.1}
  \ee

We now rewrite this in a way which is easier to deal with. Set $\mbox{ch}_{0}(\ep)=1$. Recall that
$|\bar\si|_{+}-|\bar\si| =\sum_{i=1}^{\ff}  k_{i}$. Setting 
$k_{0}=n-|\bar\si|_{+}$ we then have 
%$n-|\bar\si|=\sum_{i=1}^{\ff}  ik_{i}$ and 
$n-|\bar\si| =\sum_{i=0}^{\ff}  k_{i}$, the total number of chains when we include $k_{0}$ chains of order $0$. We can then rewrite (\ref{11.1}) as
  \begin{eqnarray}
&&
 \hspace{-.4 in}L_{1}^{n} (x,\ep)=\sum_{\{\bar\si\,|\, 0\leq |\bar\si|\leq |\bar\si|_{+} \leq n\} }{n! \over \prod_{i=0}^{\ff}k_{i}! }\prod_{i=0}^{\ff} \(\mbox{ch}_{i}( \ep)\)^{k_{i}} 
    L_{\sum_{i=0}^{\ff}  k_{i}} (x,\ep). \label{rilt.19a}
  \end{eqnarray} 

We claim that
\begin{equation}
L_{1}^{n} (x,\ep)=\sum_{k=0}^{n}{n! \over k! }\sum_{\stackrel{j_{1},\ldots,j_{k}}{\sum_{b=1}^{k}(j_{b}+1)=n } }\prod_{b=1}^{k}  \mbox{ch}_{j_{b}}( \ep)  \,\,  L_{k} (x,\ep)\label{rilt.20}
\end{equation}
where the second sum is over sequences of $k$ integers $j_{i}\geq 0$.
To see this, we simply note that if there are a total of $k$ chains of which  $k_{i}$ are of length $i$,   there are $k! /  \prod_{i=0}^{\ff}k_{i}!$ distinct ways to order them. Also, we have used the fact that
$\sum_{b=1}^{k}(j_{b}+1)=\sum_{i=1}^{\ff} (i+1)k_{i}+k_{0}=|\bar\si|_{+}+(n-|\bar\si|_{+})=n$.

A similar analysis then shows that 
\bea
&&
H_{1}^{n} (x,\ep)\label{rilt.21}\\
&&=\sum_{k=0}^{n}{n! \over k! }\sum_{r=0}^{\ff}{1 \over r! }\sum_{\stackrel{i_{1},\ldots,i_{r};\,j_{1},\ldots,j_{k}}{\sum_{a=1}^{r}i_{a}+\sum_{b=1}^{k}(j_{b}+1)=n } }
\prod_{a=1}^{r} \( {\mbox{cy}_{i_{a}}( \ep) \over 2 i_{a}}\)
\prod_{b=1}^{k}  \mbox{ch}_{j_{b}}( \ep)  \,\, \wt H_{k} (x,\ep). \nn
\eea

We can write  (\ref{rilt.21}) and (\ref{rilt.20})
as
\begin{equation}
H_{1}^{n} (x,\ep)=\sum_{k=0}^{n}A_{n,k} \wt H_{k} (x,\ep),\hspace{.2 in}L_{1}^{n} (x,\ep)=\sum_{k=0}^{n}B_{n,k} L_{k} (x,\ep).\label{15.5}
\end{equation}
where
\be
A_{n,k}  
 ={n! \over k! }\sum_{r=0}^{\ff}{1 \over r! }\sum_{\stackrel{i_{1},\ldots,i_{r};\,j_{1},\ldots,j_{k}}{\sum_{a=1}^{r}i_{a}+\sum_{b=1}^{k}(j_{b}+1)=n } }
\prod_{a=1}^{r} \( {\mbox{cy}_{i_{a}}( \ep) \over 2 i_{a}}\)
\prod_{b=1}^{k}  \mbox{ch}_{j_{b}}( \ep) \label{15.7}
\ee 
and
\begin{equation}
B_{n,k}= {n! \over k! }\sum_{\stackrel{j_{1},\ldots,j_{k}}{\sum_{b=1}^{k}(j_{b}+1)=n } }\prod_{b=1}^{k}  \mbox{ch}_{j_{b}}( \ep). \label{15.6}
\end{equation}

{\bf  Proof of Theorem \ref{theo-ILT}: } Using the above, this follows as in \cite[Section 4]{MRmem}. In somewhat more detail,  let $\MM(\HH_{K})$ denote the set of functions measurable with respect to
$\HH_{K}=: \si(H_{1}(\mu);\,\mu\in \GG_{K}^2)$. We define the ring
homomorphism \be
\Phi:\MM(\HH_{K})\mapsto \MM(\HH_{K}\times \FF)\label{mmc} \ee 
as the measurable
extension of the mapping $\Phi$ such that $\Phi(1)=1$ and
\be
\Phi\(\prod_{i=1}^n H_{1}(\mu_i)\)=\prod_{i=1}^n \(   { :G^{2}:(\mu_i)\over 2}+L_1(\mu_i) \),
\hspace{.2in}n=1,\ldots,\label{phi1}
\ee 
where $\FF$ is the $\si$-algebra generated by $X$.  With this
notation Theorem
\ref{theo-CAF} can be reformulated as follows: Let $(h_1,h_2,\ldots)$ be a sequence of
$\HH_{K}$ measurable functions. Then for any $\CC$ measurable non-negative
function $F$ on $R^{\ff}$
\be E_{G}\mathcal{P}_{\al^{2}}\(F(\Phi(h_1),\Phi(h_2),\ldots)\)
=E_{G}\(F(h_1,h_2,\ldots) \).\label{15.3}
\ee 

It follows by induction from (\ref{15.5}) that $\wt H_{n}(x,\ep)\in \HH$, and  using  (\ref{15.5}) and (\ref{phi1}) we have
\begin{equation}
\sum_{k=0}^{n}A_{n,k} \Phi(\wt H_{k} (x,\ep))=\Phi( H_{1}^{n} (x,\ep))=\({ :G^{2}:(x,\ep)\over 2}+L(x,\ep)\)^{n}.\label{15.4yu}
\end{equation}
The proof of \cite[Lemma 4.5]{MRmem} then shows that 
\begin{equation}
\Phi(\wt H_{n}(x,\ep))=  \sum_{m=0}^{n}{n \choose m}\Psi_{m}(x,\ep) L_{n-m}(x,\ep)  \label{15.4}
\end{equation}
where $\Psi_{m}(x,\ep)$ is such that  
\bea
\({ :G^{2m}:\over 2^{m}}\times L_{n-m}\)(\nu)&\stackrel{def}{=}&\lim_{\ep\to 0}\int :G^{2m}:(x,\ep) L_{n-m}(x,\ep)\,d\nu(x)\nn\\
&=&\lim_{\ep\to 0}\int \Psi_{m}(x,\ep) L_{n-m}(x,\ep)\,d\nu(x)\label{cross1}
\eea
for all $0\leq m\leq n$.
Then, using Lemma \ref{lem-comb}, the proof of our Theorem will follow from the fact that $\Phi$ is an isometry in $L^{2}$.\qed

To better appreciate the nature of $\Psi_{m}(x,\ep)$ and   (\ref{cross1}) we note that when $m=n$, (\ref{cross1}) is
\begin{equation}
 { :G^{2n}(\nu):\over 2^{n}} =\lim_{\ep\to 0}\int \Psi_{n}(x,\ep) \,d\nu(x).\label{cross2}
\end{equation}

\section{Proof of Lemma \ref{lem-comb}}\label{sec-comb}

\noindent{\bf Proof of Lemma \ref{lem-comb}:$\,$} Recall (\ref{15.5}). The key to proving our Lemma is to show that for any $\nu\in \mathcal{G}_{K}^{2n'}$
\begin{equation}
\int H_{1}^n(x,\ep)\,d\nu(x)=\sum_{k=0}^{n}A_{n,k}  H_{k} (\nu)+o(u(\ep)^{-(n'-n)}).\label{key1}
\end{equation}
Our Lemma then follows easily. See the end of the proof of \cite[Lemma 4.3]{MRmem}, starting from (4.84).

Recall that \[H_{1}(\nu)= { :G^{2}:(\nu)\over 2}+\sqrt{2}\al\,G(\nu)+\al^{2}  |\nu|.\] Therefore,
see \cite[Lemma 3.3]{MRmem}
\be H_{1}(x,\ep)\st H_{1}(f_{x,\ep}\cdot\,dx')=\lim_{\de\rar 0}\int
\({:G^2_{x',\de}:\over 2}  + 2^{1/2}\al\, G_{x',\de}+\al^{2}\)f_{x,\ep}(x')\,dx'\label{12.1} \ee where the limit is taken in $L^2$. Since it follows from \cite[Lemma 8.6]{LMR2} that 
$f_{x,\ep}\cdot\,dx'\in\GG_K^2$, $H_{1}(x,\ep)$ is one of the
basic random variables that generate $\HH$. As explained 
 in the proof of \cite[Theorem  4.2]{MRmem}, for fixed $\ep>0$, $H\st\{H_{1}(x,\ep),x\in R^m\}$
can be taken to be continuous almost surely.   Furthermore, the
convergence in (\ref{12.1}) is almost sure and in $L^p$, for all $p$, since Gaussian chaos
processes have all moments.

 Clearly 
 \be
H_{1}^n(x,\ep)=\lim_{\de\rar 0}
\prod_{i=1}^n\int \({:G^2_{x_i,\de}:\over 2}  + 2^{1/2}\al\, G_{x_i,\de}+\al^{2}\)f_{x,\ep}(x_i)\,dx_i.\label{12.1x} \ee We define
\be H^n_{1,\ep}\nu = \int H_{1}^n(x,\ep)\,d\nu(x)\label{12.3a}. \ee Since the
right-hand side of (\ref{12.1x}) converges in $L^2$ uniformly in $x$ as $\de\to 0$,
we see that
\be H^n_{1,\ep}\nu :=\lim_{\de\rar 0}\int\!\!\int \prod_{i=1}^n
\({:G^2_{x_i,\de}:\over 2}  + 2^{1/2}G_{x_i,\de}\al+\al^{2}\)f_{x,\ep}(x_i)\,dx_i\,d\nu(x)\label{12.3aa} \ee in $L^2$. (In fact by
\cite[Lemma 3.3]{Arc-Gin} it also converges almost surely and in $L^p$ for all $p\ge 0$).

Expand $\prod_{i=1}^n \({:G^2_{x_i,\de}:\over 2}  + 2^{1/2}\al\, G_{x_i,\de} +\al^{2}\)$ as a sum of Wick products. Using
\cite[(2.15)]{MRmem} we can write
 \bea &&
\prod_{i=1}^n \({:G^2_{x_i,\de}:\over 2}  + 2^{1/2}\al\, G_{x_i,\de} +\al^{2}\)\label{12.2z}\\
&&
=\sum_{A,B,C}2^{|B|/2}\al^{|B|+2|C|} \sum_{R,S,T, U,V}{1\over 2^{|R|}}\nn\\
&& 
\sum_{\stackrel{\stackrel{pairings\,\PP}{of\, (R\times\{1,2\})
\cup S\cup U}} {\wt{\PP}_{k,1}\neq
\wt{\PP}_{k,2}}}\prod_{k=1}^{|R|+(|S|+|U|)/2}u_{\de,\de}(x_{\wt{\PP}_{k,1}}
-x_{\wt{\PP}_{k,2}})
 :\prod_{i\in T} {G^2_{x_i,\de}\over 2}
\prod_{j\in S\cup V}G_{x_j,\de}: \nn
\eea 
where the first sum runs over all partitions $A\cup B\cup C=\{1,2,\ldots,n\}$, and the second sum runs over all partitions $R\cup S\cup T=A$, $U\cup V=B$
  with $|S|+|U|$ even , and
the third sum runs over all pairings $\PP$ of the set $(R\times\{1,2\})\cup S\cup U$ such
that
$\wt{\PP}_{k,1}\neq \wt{\PP}_{k,2}$, where letting $(\PP_{k,1},\PP_{k,2})$ denote
the $k-$th pair of the pairing $\PP$, we set $\wt{\PP}_{k,1}=i$ if either
$\PP_{k,1}=i\times 1$ or $i\times 2$ for
$i\in R$, or 
$\PP_{k,1}=i$ for $i\in S \mbox{ or }U$, and similarly for $\wt{\PP}_{k,2}$.  Here
we use the fact that for $i\in S$ one of the two $G_{x_i,\de}$ terms is allocated to
the $u_{\de,\de}$ terms and the other to the Wick product. Since there are
two ways to do this the 1/2 is cancelled.  (\ref{12.2z}) can be rewritten as
\bea
&&\hspace{-.5 in}
\prod_{i=1}^n \({:G^2_{x_i,\de}:\over 2}  + 2^{1/2}\al\, G_{x_i,\de} +\al^{2}\)=\sum_{R,S,T, U,V,C} 2^{(|U|+|V|)/2}\al^{|U|+|V|+2|C|}\nn \\
&&\hspace{.8 in}
\EE_{\de}(x_1,\ldots,x_n;R,S, U) :\prod_{i\in T} {G^2_{x_i,\de}\over 2}
\prod_{j\in S\cup V}G_{x_j,\de}:\label{12.2}
\eea where  \be
\EE_{\de}(x_1,\ldots,x_n;R,S, U)
\st{1\over 2^{|R|}}\sum_{\stackrel{\stackrel{pairings\,\PP}{of\, (R\times\{1,2\})
\cup S\cup U}} {\wt{\PP}_{k,1}\neq
\wt{\PP}_{k,2}}}\prod_{k=1}^{|R|+(|S|+|U|)/2}u_{\de,\de}(x_{\wt{\PP}_{k,1}}
-x_{\wt{\PP}_{k,2}}).
\ee

Using this, (\ref{12.3aa}), Fubini's theorem and the definition of $f_{x,\ep}$ we
have
\bea 
H^n_{1,\ep}\nu  &=&\sum_{R,S,T, U,V,C} 2^{(|U|+|V|)/2}\al^{|U|+|V|+2|C|}
\lim_{\de\rar 0}\int\EE_{\de}(x_1,\ldots,x_n;R,S, U)
\nn\\ &&\(\int :\prod_{i\in T}
{G^2_{x+x_i,\de}\over 2}
\prod_{j\in S\cup V}G_{x+x_j,\de}:\,d\nu(x)\)\prod_{i=1}^n f_{\ep}(x_i)\,dx_i.\label{12.3}
\eea

Let us now consider
\bea &&
\int \EE_{\de}(x_1,\ldots,x_n;R,S,U)
\prod_{i=1}^n f_{\ep}(x_i)\,dx_i\label{12.6}\label{zz}\\ && \qquad={1\over
2^{|R|}}\sum_{\stackrel{\stackrel{pairings\,\PP}{of\,(R\times\{1,2\})\cup S\cup U}}
{\wt{\PP}_{k,1}\neq
\wt{\PP}_{k,2}}}\int
\prod_{k=1}^{|R|+(|S|+|U|)/2}u_{\de,\de}(x_{\wt{\PP}_{k,1}}-x_{\wt{\PP}_{k,2}})
\prod_{i=1}^n f_{\ep}(x_i)\,dx_i.\nn
\eea 
%Suppose $|R|+(|S|+|U|)/2=q$. Then by the multiple H\"older inequality and
%\cite[(3.7)]{MRmem}, (\ref{zz}) is bounded above by
% \bea &&C\int\int
%|u_{\de,\de'}(|x-y|)|^q f_{\ep}(x)f_{\ep}(y)\,dx\,dy\nn\\ &&\qquad\le C\int\int
%|u(|x-y+(x'-y')|)|^q f_{\ep}(x)f_{\ep}(y)\,dx\,dy\rho_\de(dx')
%\rho_\de(dy').\nn
%\eea It follows from this and  \cite[Lemma 3.2]{MRmem} that (\ref{zz}) is
%$O((u(\ep))^{(|R|+(|S|+|U|)/2)})$ for all $\de>0$. 
We   reorganize this in a form which is more useful. Fix some pairing $\PP$ in
the sum and pick any factor $u_{\de,\de}(x_i-x_j)$ in the product corresponding to
$\PP$. If both $i,j\in S\cup U$ we think of
$i,j$ as forming a chain of order one. $u_{\de,\de}(x_i-x_j)$ is the factor associated with
this chain. If say $j\in R$, there will be one other factor in the product
corresponding to
$\PP$ which contains $x_j$, say $u_{\de,\de}(x_j-x_k)$. If both $i,k\in S\cup U$, we think of $i,j,k$
as forming a chain of order two. $u_{\de,\de}(x_i-x_j) u_{\de,\de}(x_j-x_k)$ is the factor
associated with this chain. If either $i$ or $k$ or both are in $R$, we continue to
find the other factors containing them, and continue in this manner until we  can
go no further. Two possibilities arise. Either we end up with a chain of elements
$i_1,i_2,\ldots,i_v$ with end points $i_1,i_v\in S\cup U$ and intermediate points
$i_2,\ldots,i_{v-1} \in R$ and associated factor
\be
\prod_{j=2}^v u_{\de,\de}(x_{i_{j}}-x_{i_{j-1}})\label{12.7} \ee (such a chain is said to be of
order $v-1$ ), or we have, what we call,  a cycle
$i_1,i_2,\ldots,i_v$ with all elements in $R$  and associated factor
\be u_{\de,\de}(x_{i_{1}}-x_{i_{v}})\prod_{j=2}^v u_{\de,\de}(x_{i_{j}}-x_{i_{j-1}}) \label{12.8}
\ee (such a cycle is  said to be of order $v$).

In this way the product
\be
\prod_{k=1}^{|R|+(|S|+|U|)/2}u_{\de,\de}(x_{\wt{\PP}_{k,1}}-x_{\wt{\PP}_{k,2}}) \ee in
(\ref{12.4}) associated with $\PP$ breaks up into a product of factors associated
with the chains and cycles of $\PP$. Note that each $\PP$ appearing in (\ref{12.4})
will necessarily have precisely $(|S|+|U|)/2$ chains. Set $p=(|S|+|U|)/2$.

When
$\PP$ decomposes into $k_l$ chains of order
$l,\,l=1,2,\ldots$ and 
$m_l$ cycles of order $l,\,l=2,\ldots$, we write $\PP\rar \si=(k_1,\ldots;\,m_{2},\ldots)$. Letting $\mbox{ch}_{j,\de}(\ep)$ and
$\mbox{cy}_{j,\de}(\ep)$ denote  the chain and cycle factors
defined in (\ref{1L.1}), (\ref{4.2}) with $u$ replaced by $u_{\de,\de}$, we have
\bea && {1\over 2^{|R|}}\sum_{\stackrel{\stackrel{pairings\,\PP}{of\, (R\times\{1,2\})\cup S\cup U}}
{\wt{\PP}_{k,1}\neq
\wt{\PP}_{k,2}}}\int
\prod_{k=1}^{|R|+p}u_{\de,\de}(x_{\wt{\PP}_{k,1}}-x_{\wt{\PP}_{k,2}}) 
\prod_{i\in R\cup S\cup U}f_{\ep}(x_i)\,dx_i \nn\\ && ={1\over
2^{|R|}}\sum_{\si=(k_1,\ldots;\,m_{2},\ldots)}\sum_{\stackrel{\stackrel
{pairings\,\PP}{of\,(R\times\{1,2\})\cup S\cup U}} {\PP\rar
\si}}\prod_{l=1}^\ff\,
(\mbox{ch}_{l,\de}(\ep))^{k_l}\,(\mbox{cy}_{l,\de}(\ep))^{m_l}\label{12.10}.
\eea
There are no cycles of order one. The notation for the last product is purely for convenience, and we take $\mbox{cy}_{1,\de}(\ep)=1, m_{1}=0$.

Note that when $\PP\rar \si=(k_1,\ldots;\,m_{2},\ldots)$ and $\PP$ is a pairing of
$(R\times\{1,2\})\cup S\cup U$ we must have
$|\si|_{+}=\sum_{l=1}^\ff k_l(l+1)+\sum_{l=2}^\ff m_ll=|R|+|S|+|U|$. We now further simplify
(\ref{12.10}) by observing that the number of pairings $\PP$ of 
$(R\times\{1,2\})\cup S\cup U$ with $\PP\rar \si=(k_1,\ldots;\,m_{2},\ldots)$ is
 \bea &&
\label{12.11}\\ && {|R|!\over \prod_{l=2}^\ff (l!)^{m_l}(m_l!)((l-1)!)^{k_l}
(k_l!) }\prod_{l=1}^\ff ({(l-1)!\over 2})^{m_l} ((l-1)!)^{k_l}{(|S|+|U|)!\over 
2^{\sum_{l=2}^\ff
k_l}}2^{|R|}.\nn
\eea 
Here, the first factor gives the number of ways to partition $R$ into
$m_l$ cycles of length $l,\,l=2,\ldots$ and $k_l$ mid-chains (i.e. chains with
end points deleted) of length
$l-1,\,l=1,\ldots$. To get the remainder of (\ref{12.11}) we note that in each cycle
of length $l$ we can permute the points of the cycle in
$(l-1)!$ distinct ways, except that we must divide by $2$ to take into  account the
mirror image if $l>2$, (we will explain shortly where the factor $1/2$ for cycles of
length
$l=2$ comes from), while for each chain of length $l$ we can permute the elements
of the mid-chain in $(l-1)!$ ways, and the $|S|+|U|$ end points can be permuted among
themselves in
$(|S|+|U|)!$ ways, except that for any of the $p=\sum_{l=1}^\ff k_l$ given chains we
mustn't count an interchange of the end points of the same chain, since  that has
already been counted when we considered the permutations of the mid-chain.
Finally,  recall that the pairings are actually parings of
$(R\times\{1,2\})\cup S\cup U$, not of $R\cup S\cup U$, so that for any given  pairing we can
get analogous but distinct pairings by interchanging $i\times 1$ with 
$i\times 2$ for each $i\in R$. The only exception is that for any cycle of length $l=2$
we get $2$ rather than $4$ distinct pairings. Altogether this gives rise to $2^{|R|}/
2^{m_2} $ distinct pairings. (This explains where  the factor $1/2$ for cycles of
length
$l=2$ in (\ref{12.11}) comes from). Therefore, combining (\ref{12.10}) and
(\ref{12.11}) we see that
\bea 
&&\int \EE_{\de}(x_1,\ldots,x_n;R,S,U)
\prod_{i=1}^n f_{\ep}(x_i)\,dx_i\label{12.12}\\&&= {1\over 2^{|R|}}\sum_{\stackrel{\stackrel{pairings\,\PP}{of\, (R\times\{1,2\})\cup S\cup U}}
{\wt{\PP}_{k,1}\neq
\wt{\PP}_{k,2}}}\int
\prod_{k=1}^{|R|+p}u_{\de,\de}(x_{\wt{\PP}_{k,1}}-x_{\wt{\PP}_{k,2}}) 
\prod_{i\in R\cup S\cup U}f_{\ep}(x_i)\,dx_i \nn\\
 && =\sum_{\stackrel{\si=(k_1,\ldots;\,m_{2},\ldots)}{|\si|_{+}=|R|+|S|+|U|}}{2^{ -\sum_{l=1}^\ff
k_l}|R|!\over \prod_{l=1}^\ff
(l!)^{m_l}(m_l!)((l-1)!)^{k_l}(k_l!) }
\nn\\ &&\qquad\qquad
\prod_{l=1}^\ff ({(l-1)!\over 2})^{m_l} ((l-1)!)^{k_l}(|S|+|U|)!\prod_{l=1}^\ff \,(ch^{\ep}_{l,\de})^{k_l}\,
(\mbox{cy}_{l,\de}(\ep))^{m_l} \nn\\ &&  
\hspace{-.3 in}=\sum_{\stackrel{\si=(k_1,\ldots;\,m_{2},\ldots)}{\underset{2p=|S|+|U|}{|\si |_{+}=|R|+|S|+|U|}}}{2^{-(|S|+|U|)/2}|R|!(|S|+|U|)!\over \prod_{l=1}^\ff (m_l!) (k_l!) }
\prod_{l=1}^\ff \,(\mbox{ch}_{l,\de}(\ep))^{k_l}\,\({\mbox{cy}_{l,\de}(\ep) \over 2l}\)^{m_l}.\nn
\eea

Let 
    \begin{equation}
 h(s)   :=  \int_{|\xi|\leq s}  1/\psi(|\xi|)  \,d\xi.\label{dhr.1}
   \end{equation}
   By (\ref{sl1.5}), $\lim_{s\to\ff}h( s)  =\ff.$ We mention some results which are analogues of results used in \cite{MRmem}. We will prove these results under the assumptions of this paper at the end of this section.
   
    \bl\label{lem-chat} For any $\de>0$
  \begin{equation}
\mbox{ch}_{k,\de}(\ep)\leq \mbox{ch}_{k}(\ep)= 
 O\(   (h(1/\ep))^{k}\)   \qquad\mbox{as }\ep\rar 0\label{ss}
  \end{equation}
   and
   \begin{equation}
\mbox{cy}_{k,\de}(\ep)= \mbox{cy}_{k}(\ep)= 
O\(   (h(1/\ep))^{k}\)   \qquad\mbox{as }\ep\rar 0.\label{ssci}
  \end{equation}
  In addition, for $\ep>0$ fixed
  \begin{equation}
 \lim_{\de\to 0}\mbox{ch}_{k,\de}(\ep)=\mbox{ch}_{k}(\ep) \hspace{.2 in}\mbox{and}\hspace{.2 in} 
  \lim_{\de\to 0}\mbox{cy}_{k,\de}(\ep)=\mbox{cy}_{k}(\ep).\label{ssd}
  \end{equation}
 \el
 
   \bl\label{lem-ghat} For $\nu\in\GG_{K}^{2n'}$
   \begin{equation}
  \sup_{|x_i|\le \ep}\| \int :\prod_{i=1}^{ k}G_{x+x_i,\de}:\,d\nu(x)-:G^{k}_{\de}\nu:\|_2 
  =o\left((h(1/\ep))^{-( n'-k/2)}   \right)   \qquad\mbox{as }\ep\rar 0.\label{sse}
   \end{equation}
   \el
   
   Since $\sum_{l=1}^\ff k_ll+\sum_{l=2}^\ff m_ll=|R|+(|S|+|U|)/2$ it follows from (\ref{12.12}) and Lemma \ref{lem-chat} that 
   \begin{equation}
   \int \EE_{\de}(x_1,\ldots,x_n;R,S, U) \prod_{i=1}^n
f_{\ep}(x_i)\,dx_i=O\(   (h(1/\ep))^{|R|+(|S|+|U|)/2}\) \label{ssf}
   \end{equation}
as $\ep\rar 0$.

Using this and  Lemma \ref{lem-ghat}  we see that for $\nu\in\GG_{K}^{2n'}$
\bea &&
\| \int \EE_{\de}(x_1,\ldots,x_n;R,S,U)\label{tot.1}\\ &&\quad
\lc \int :\prod_{i\in T} {G^2_{x+x_i,\de}\over 2}\prod_{j\in
S\cup V}G_{x+x_j,\de}:\,d\nu(x)-{:G^{2|T|+|S|+|V|}_{\de}\nu:\over 2^{|T|}}\rc\prod_{i=1}^n f_{\ep}(x_i)\,dx_i\|_2\nn\\ &&
\leq \int\EE_{\de}(x_1,\ldots,x_n;R,S,U) \prod_{i=1}^n f_{\ep}(x_i)\,dx_i \nn\\
&&\quad
\sup_{|x_i|\le \ep}\| \int :\prod_{i\in T} {G^2_{x+x_i,\de}\over 2}\prod_{j\in
S\cup V}G_{x+x_j,\de}:\,d\nu(x)-{:G^{2|T|+|S|+|V|}_{\de}\nu:\over 2^{|T|}}\|_2\nn\\ &&
\leq  O\(   (h(1/\ep))^{|R|+(|S|+|U|)/2}\)o(( (h(1/\ep))^{-(n'-(|T|+(|S|+|V|)/2))}\nn\\ && =o( (h(1/\ep))^{-(n'-n)})\nn
\eea 
for all $\de>0$. 

Note that it follows from (\ref{12.12}) and Lemma \ref{lem-chat} that
 \bea &&
\lim_{\de\to 0}\int \EE_{\de}(x_1,\ldots,x_n;R,S, U) \prod_{i=1}^n f_{\ep}(x_i)\,dx_i\label{zz1}\\ &&\quad =\sum_{\stackrel{\si=(k_1,\ldots;\,m_{2},\ldots)}{\underset{2p=|S|+|U|}{|\si |_{+}=|R|+|S|+|U|}}}{2^{-(|S|+|U|)/2}|R|!(|S|+|U|)!\over \prod_{l=1}^\ff (m_l!) (k_l!) }
\prod_{l=1}^\ff \, (\mbox{ch}_{l}(\ep))^{k_l}\,\({\mbox{cy}_{l}(\ep) \over 2l}\)^{m_l}.\nn
\eea 
Using this, (\ref{tot.1}) and (\ref{12.3}) we see that for $\nu\in
\GG_{K}^{2\si}$
\bea 
&&\hspace{-.2 in}
H^n_{1,\ep}\nu =
\sum_{R,S,T, U,V, C}  2^{(|U|+|V|)/2}\al^{|U|+|V|+2|C|}  \lim_{\de\rar 0}\int\EE_{\de}(x_1,\ldots,x_n;R,S,U)\nn\\
&&\prod_{i\in R\cup S\cup U}f_{\ep}(x_i)\,dx_i :{G^{2|T|+|S|+|V|}_{\de}\over
2^{|T|}}\nu:+o((u(\ep))^{-(n'-n)}) \label{12.4}\\ 
&&\hspace{-.2 in}=\sum_{k=0}^{2n}\sum_{\stackrel{R,S,T, U,V, C}{2|T|+|S|+|V|=k}}2^{(|U|+|V|)/2}\al^{|U|+|V|+2|C|}  \lim_{\de\rar 0}
\int\EE_{\de}(x_1,\ldots,x_n;R,S,U)\nn\\ &&\qquad\prod_{i\in R\cup
S\cup U}f_{\ep}(x_i)\,dx_i {1\over 2^{|T|}}:G^{k}_{\de}\nu:+o((u(\ep))^{-(n'-n)})\nn\\
&=&\sum_{k=0}^{2n}\sum_{\stackrel{R,S,T, U,V, C}{2|T|+|S|+|V|=k}}2^{(|U|+|V|)/2}\al^{|U|+|V|+2|C|}  \nn\\
&&\hspace{-.2 in}\sum_{\stackrel{\si=(k_1,\ldots;\,m_{2},\ldots)}{\underset{2p=|S|+|U|}{|\si |_{+}=|R|+|S|+|U|}}}{2^{-(|S|+|U|)/2}|R|!(|S|+|U|)!\over \prod_{l=1}^\ff (m_l!) (k_l!) }
\prod_{l=1}^\ff \, (\mbox{ch}_{l}(\ep))^{k_l}\,\({\mbox{cy}_{l}(\ep) \over 2l}\)^{m_l}{:G^{k}\nu:\over
2^{|T|}}\nn\\ &&\hspace{3.5 in}+o((u(\ep))^{-(n'-n)})\nn \eea 
in $L^2$, as $\ep\to 0$. We remark that 
if $|R|+|S|+|U|=0$ we will have $|\si|=0$.

Recalling the definition of $I_{n,\ep} (\si)$, see (\ref{4.5}),
and then combining (\ref{12.12}) with (\ref{12.4}) we see that for $\nu\in
\GG_{K}^{2n'}$
\bea 
H^n_{1,\ep}\nu 
&=&\sum_{k=0}^{2n}\sum_{\stackrel{R,S,T, U,V, C}{2|T|+|S|+|V|=k}}2^{(|U|+|V|)/2-(|S|+|U|)/2}\al^{|U|+|V|+2|C|}   {|R|!(|S|+|U|)! \over n!} \nn\\  &&\hspace{-.7 in}
\sum_{\stackrel{\{\si\,|\, 0\leq |\si| \leq |\si|_{+} \leq n\} }{|\si|_{+}=|R|+|S|+|U|,\,2p=|S|+|U|}} (n-|\si|_{+})!I_{n,\ep} (\si){:G^{k}\nu:\over
2^{|T|}}+o((u(\ep))^{-(n'-n)}).\label{12.4y} \eea 
in $L^2$, as $\ep\to 0$.

Let us introduce the abbreviation
\begin{equation}
Z( k, |\si|_{+},p)=\{   2|T|+|S|+|V|=k,\,|R|+|S|+|U|=|\si|_{+},\,|S|+|U|=2p\}.\label{}
\end{equation}
Since
\begin{eqnarray}
&&(|U|+|V|)/2-(|S|+|U|)/2-|T|
\label{k1}\\
&&=|V|-(2|T|+|S|+|V|)/2=|V|-k/2   \nonumber
\end{eqnarray} 
 and 
 \begin{eqnarray}
 &&|U|+|V|+2|C|
 \label{k2}\\
 &&=2n-2|R|-2|S|-2|T|-|U|-|V|\nn\\
 &&=2n-2\(|R|+|S|+|U|\)+\(|S|+|U|\)-\(2|T|+|S|+|V|\)   \nonumber\\
 &&=2n-2|\si|_{+}+2p-k=2(n-|\si|)-k,   \nonumber
 \end{eqnarray}
we can reorganize (\ref{12.4y}) as 
\bea 
H^n_{1,\ep}\nu 
&=&\sum_{\{\si\,|\, 0\leq |\si| \leq |\si|_{+} \leq n\} }\,\,\sum_{k=0}^{2n}\sum_{\stackrel{R,S,T, U,V, C}{Z( k, |\si|_{+},p)}}\label{12.4yy1}\\  && \hspace{-.7 in}\al^{2(n-|\si|)-k }  2^{|V|} {|R|!(|S|+|U|)! \over n!}{:G^{k}\nu:\over
2^{k/2}} (n-|\si|_{+})!I_{n,\ep} (\si)+o((u(\ep))^{-(n'-n)}).\nn
\eea
Since there are $n!/(|R|!|S|!|T|!|U|!|V|! |C|!)$  ways to partition $[1,n]$ into sets of size $|R|,|S|,|T|, |U|,|V|, |C|$ we see that
\bea
H^n_{1,\ep}\nu 
&=&\sum_{\{\si\,|\, 0\leq |\si| \leq |\si|_{+} \leq n\} }\,\,\sum_{k=0}^{2n}\sum_{\stackrel{|R|,|S|,|T|, |U|,|V|, |C|}{Z( k, |\si|_{+},p)}}{(n-|\si|_{+})! \over |T|!  |V|!  |C|! }\label{12.4yy}\\  && \al^{2(n-|\si|)-k }  2^{|V|} {|S|+|U| \choose |S|} {:G^{k}\nu:\over
2^{k/2}}I_{n,\ep} (\si)+o((u(\ep))^{-(n'-n)}).\nn
 \eea 
Now  write  (\ref{12.4yy}) as 
\bea
&&
\hspace{-.2 in}H^n_{1,\ep}\nu =\sum_{\{\si\,|\, 0\leq |\si| \leq |\si|_{+} \leq n\} }\(\sum_{k=0}^{2n}\rho(\si,k)\al^{2(n-|\si|)-k }  {:G^{k}\nu:\over
2^{k/2}}\)I_{n,\ep} (\si)\nn\\
&&\hspace{ 2.5in}+o((u(\ep))^{-(n'-n)}),\label{12.4w}
\eea
where
\begin{equation}
\rho(\si,k)=\sum_{\stackrel{|R|,|S|,|T|, |U|,|V|, |C|}{Z( k, |\si|_{+},p)}}{ (n-|\si|_{+})!\over |T|!  |V|!  |C|! }  {|S|+|U| \choose |S|}2^{|V|}.\label{12.4x}
\end{equation}

\bl\label{lem-comba}
\begin{equation}
\rho(\si,k)={2(n-|\si|) \choose k}.\label{12.4z}
\end{equation}
\el

Using this in (\ref{12.4w}) and recalling the definition (\ref{4.0}) of $H_{n-|\si|}(\nu)$ we have 
\begin{equation}
H^n_{1,\ep}\nu =\sum_{\{\si\,|\, 0\leq |\si| \leq |\si|_{+} \leq n\} }I_{n,\ep} (\si)H_{n-|\si|}(\nu)+o((u(\ep))^{-(n'-n)}).\label{12.4g}
\end{equation}
  Then by the argument which led to (\ref{15.5}) we obtain (\ref{key1}). As explained in the beginning of this section, (\ref{key1}) will complete the proof of Lemma \ref{lem-comb}.

{\bf  Proof of Lemma \ref{lem-comba}: }
Under our constraint $Z( k, |\si|_{+},p)$
\[2|T|+|S|+|V|=k,\,|R|+|S|+|U|=|\si|_{+},\,|S|+|U|=2p.\] 
Once we have specified $0\leq |S|\leq 2p$ and $0\leq |V|\leq k-|S|$, then $|U|,|R|,|T|$ are determined and consequently so is $|C|$. Furthermore, we have
\begin{equation}
|T|+|V|+|C|=n-(|R|+|S|+|U|)=n-|\si|_{+},\label{}
\end{equation}
$|T|=(k-|S|-|V|)/2$ and 
\bea
&&
|C|=n-(|R|+|S|+|T|+ |U|+|V|)\label{ccc}\\
&&=n-|\si|_{+}- (|T|+|V|)= n-|\si|_{+}-(k-|S|+|V|)/2. \nn
\eea
 Thus, setting $m=n-|\si|_{+}$ we have
 \begin{eqnarray}
 &&\hspace{ -.5in}\rho(\si,k)=\sum_{|S|=0}^{2p}{2p \choose |S|}
 \label{12.4h}\\
 &&\hspace{ .4in}\sum_{|V|=0}^{k-|S|}{m \choose (k-|S|-|V|)/2,\hspace{.1 in}|V|,\hspace{.1 in} m-(k-|S|+|V|)/2}2^{|V|}   \nonumber
 \end{eqnarray}
and our Lemma then follows if we can show  that for all  $k,m,p\in Z_{+}$ 
\begin{equation}
\sum_{s=0}^{2p}{ 2p\choose s}\sum_{v=0}^{k-s}{ m\choose (k-s-v)/2,\hspace{.1 in} v,\hspace{.1 in} m-(k-s+v)/2}2^{v}=  {2m+2p \choose k},\label{12.4j}
\end{equation}
where the sum over $v$ is only taken over those values of $v$ such that $(k-s-v)/2$ and $m-(k-s+v)/2$ are non-negative integers.

To see (\ref{12.4j}), note  first that  
  \[
   \binom{2m+2p}{k} = \sum_{s=0}^{2p}\binom{2p}{s} \binom{2m}{k-s},
  \]
which can be seen by examining the coefficient of $x^{k}$ on both sides of
  \[
   (x+1)^{2m+2p} = (x+1)^{2p}(x+1)^{2m}.
  \]
Thus, it will suffice to prove that
  \[
   \sum_{v=0}^{k-s} \binom{m}{(k-s-v)/2,\hspace{.1 in}v,\hspace{.1 in}m-(k-s+v)/2} 2^v = \binom{2m}{k-s},
  \]
which we rewrite with $q=k-s$ as
  \[
   \sum_{v=0}^{q} \binom{m}{(q-v)/2,\hspace{.1 in}v,\hspace{.1 in}m-(q+v)/2} 2^v = \binom{2m}{q},
  \]
Now, the coefficient of $x^q$ in $(x+1)^{2m}$ is $\binom{2m}{q}$. But also, as
  \begin{multline*}
   (x+1)^{2m} = (x^2+2x+1)^m 
	= \sum_{i,v} \binom{m}{i,\hspace{.1 in}v,\hspace{.1 in}m-i-v} (x^2)^i (2x)^v (1)^{m-i-v}\\
= \sum_{i,v} \binom{m}{i,\hspace{.1 in}v,\hspace{.1 in}m-i-v} 2^v  x^{2i+v},
  \end{multline*}
we see that the only $(v,i)$ that contribute to $x^q$ are those with $0\leq v \leq q$ and $i=(q-v)/2$, and then only if $q$ and $v$ have the same parity. Thus, the coefficient of $x^q$ is simply
  \bea
  &&
   \sum_{v=0}^q \binom{m}{(q-v)/2,\hspace{.1 in}v,\hspace{.1 in}m-(q-v)/2-v} 2^v \label{bry5}\\
   &&\hspace{1 in}=   \sum_{v=0}^{q} \binom{m}{(q-v)/2,\hspace{.1 in}v,\hspace{.1 in}m-(q+v)/2} 2^v  .\nn
  \eea

This completes the proof of Lemma \ref{lem-comba} and hence of  Lemma \ref{lem-comb}.
\qed

{\bf  Proof of Lemma \ref{lem-chat}: } The bounds on $\mbox{ch}_{k}(\ep),\mbox{cy}_{k}(\ep)$ come from   \cite[Lemma 8.6]{LMR2}. These bounds are obtained by rewriting $\mbox{ch}_{k}(\ep),\mbox{cy}_{k}(\ep)$ in terms of the Fourier transform of $u$. Since $u_{\de,\de}=f_{\de}\ast u\ast f_{\de}$, we have
\begin{equation}
\hat u_{\de,\de}( \la)=(\hat f(\de \la))^{ 2}\,  \hat u( \la).\label{hat.1}
\end{equation}
Since by assumption $\int f( x)\,dx=1$, we have $|\hat f(\de \la)|\leq 1$ for all $\de, \la$ so that
\begin{equation}
|\hat u_{\de,\de}( \la)|\leq  \hat u( \la);\label{hat.2}
\end{equation}
 and $\lim_{\de\rar 0}\hat f(\de \la)=\hat f(0)=1$ for each $\la$. Using this in the proof of \cite[Lemma 8.6]{LMR2} the rest of our Lemma follows easily.\qed

{\bf  Proof of Lemma \ref{lem-ghat}: }Since $:G^{k}_{\de}\nu:=\int :G^{ k}_{x,\de}:\,d\nu(x)$ we can write
\bea
&&
\int :\prod_{i=1}^{ k}G_{x+z_i,\de}:\,d\nu(x)-:G^{k}_{\de}\nu: \label{ghat.10}\\
&&=\sum_{j=1}^{ k}\int :\left(   \prod_{i=1}^{ j-1}G_{x+z_i,\de}\right)(G_{x+z_j,\de}-G_{x,\de})G^{ k-j}_{x,\de}:\,d\nu(x)\nn\\
&&=\sum_{j=1}^{ k}\int :\left(   \prod_{i=1}^{ j-1}G_{x+z_i,\de}\right)(\De_{j,x}G_{x,\de})G^{ k-j}_{x,\de}:\,d\nu(x)\nn
\eea
where $\De_{j,x}f( x)=f( x+z_j)-f( x)$. We can then compute
\be 
E\left(\left( \int :\left(   \prod_{i=1}^{ j-1}G_{x+z_i,\de}\right)(\De_{j}G_{x,\de})G^{ k-j}_{x,\de}:\,d\nu(x) \right) ^{ 2}  \right)
\label{ghat.11}
\ee
as a sum of terms of the form
 \be
\int\left( \De_{j,x}\De_{j,y}u_{\de,\de}( (x-y )+a_1   )\,\,\prod_{l=2}^{ k}u_{\de,\de}( (x-y )+a_l  )\right)\,d\nu ( x)\,d\nu ( y)\label{ghat.12}
\ee
or
 \be
\hspace{-.2 in}\int\left( \De_{j,x}u_{\de,\de}( (x-y )+a_1   )\De_{j,y}u_{\de,\de}( (x-y )+a_2   )\,\,\prod_{l=3}^{ k}u_{\de,\de}( (x-y )+a_l   )\right)\,d\nu ( x)\,d\nu ( y),\label{ghat.12a}
\ee
where the $a_l$ can be either $0$ or some combination   of the $z_{m}, m=1,\ldots, k$. 
We consider (\ref{ghat.12a}). (\ref{ghat.12}) is similar and easier. Write
\begin{equation}
u_{\de,\de}( (x-y )+a_l   )=\int e^{ i\la_{l}( (x-y )+a_l )}\hat u_{\de,\de}(\la_{l} )\,d\la_{l},\label{ghat.15}
\end{equation}
\begin{equation}
 \De_{j,x}u_{\de,\de}( (x-y )+a_1   )=\int \left( e^{ i\la_{1}x_{j}} -1 \right)e^{ i\la_{1}( (x-y )+a_1)}\hat u_{\de,\de}(\la_{1} )\,d\la_{1},\label{ghat.15a}
\end{equation}
and similarly for $\De_{j,y}u_{\de,\de}( (x-y )+a_2   )$. Multiplying these out as in (\ref{ghat.12a}) and then integrating with respect to $d\nu ( x)\,d\nu ( y)$ we see that (\ref{ghat.12a}) is bounded by
\begin{equation}
\int |e^{ i\la_{1}x_{j}} -1|\,|e^{ i\la_{2}x_{j}} -1|\,\, |\hat \nu\( \la_{1}+\cdots+  \la_{k}\) |^{ 2}\,\,\prod_{l=1}^{ k}\hat u_{\de,\de}(\la_{l} )\,d\la_{l}.\label{ghat.16}
\end{equation}
Using (\ref{hat.2}) and the Cauchy-Schwarz inequality, it suffices to bound
\bea
&&
\int |e^{ i\la_{1}x_{j}} -1|^{ 2}\,\,|\hat \nu\( \la_{1}+\cdots+  \la_{k}\) |^{ 2}\,\,\prod_{l=1}^{ k}\hat u(\la_{l} )\,d\la_{l}\label{ghat.17}\\
&&=\int |\hat \nu (\la)|^{2}\,\,\tau_ {k ,r'z'}(\la)\,d\la, \nn
\eea
where
\be 
\tau_{k,r'z'}(\la)= \int |e^{ i\la_{1}x_{j}} -1|^{ 2}\, |\hat u (\la_{1}) |\tau_{ k-1}(\la-\la_{1})\,d\la_{1},
\label{2.65}
\ee 
and $\tau_{ k-1}$ is the $k-1$ fold convolution of $\hat u$.
  Our Lemma then follows from \cite[Lemma 8.4]{LMR2}.\qed

  \section{Decomposition of intersection local times}\label{sec-dec}

 The     
 renormalized intersection local times $L_{n}( \nu)$, 
   a renormalized   limit of    
  \be
  \int  \(\sum_{\om\in\mathcal{I}_{\al}} \int f_{\ep}(Y_{t}(\om)- x)\,dt \)^{n}
\,d\nu(x),\label{21.0}
\ee 
 involves both self-intersections of paths in the Poisson process $\mathcal{I}_{\al}$ and  intersections between different paths in $\mathcal{I}_{\al}$. In this section we show how to make this explicit.

 We first mimic the construction of $L_{n}( \nu)$, but for a single path under the quasi-process measure $\mu_{m}$. Let 
\begin{equation}
\LL_{1}(x,\ep):=\int  f_{\ep,x} (Y_{t})\,dt. \label{21.1}
\end{equation} 
Using the notation from Section \ref{sec-rilt} we   define  recursively 
  \begin{eqnarray}
&&
\LL_{n}(x,r)  = \LL_{1}^{n} (x,r)  -\sum_{\{ \si\,|\, 1\leq | \si| < | \si|_{+} \leq n\} } I_{n}( \si,r) \label{21.2},
  \end{eqnarray}
  where
\begin{equation}
I_{n}( \si,r)={n! \over \prod_{i=1}^{\ff}k_{i}! (n-| \si|_{+})!}\prod_{i=1}^{\ff} \(\mbox{ch}_{i}( r)\)^{k_{i} } 
    \LL_{n-|\si|} (x,r) .\label{21.3}
\end{equation} 
(Note that $n-|\si|_{+} \ge 0$.)  

Set  
\be \LL_{n,r}(\nu)=  \int   \LL_{n}(x,r) \,d\nu(x).\label{21.4}
\ee
Following the  proof  of Theorems   \ref{theo-multiriltintro-m} and  \ref{theo-1.3} we can show that if $ \nu\in \mathcal{G}_{K}^{2n}$ then 
 \begin{equation}
\LL_{n}(\nu):=\lim_{\ep \rar 0}\LL_{n,\ep}(\nu)\quad\mbox{    exists in all  $L^{p}(\mu_{m} )$},\label{21.4}
\end{equation}
and   if $n=n_{1}+\cdots+n_{k}$     and $ \nu_{i}\in \mathcal{G}_{K}^{2n_{i}}$, then 
\begin{equation}
\mu_{m}\(  \prod_{i=1}^{k} \LL_{n_{i}}(\nu_{i} ) \)= \prod_{i=1}^{k}n_{i}!   \int  \sum_{\pi\in \mathcal{M}_{a}}  \prod_{i=1}^{n-1}u (x_{\pi(i )},x_{\pi(i+1)}) \prod_{m=1}^{k}\,d\nu_{m}  (x_{m} )\label{21.4b}
\end{equation}
 where  $\mathcal{M}_{a}$ is the set of maps $\pi:[1,n]\mapsto [1,k] $  with
 $|\pi^{-1}(m)|=n_{m}$ for each $m$  and  such that    for each $l$,  if $\pi(i)=m$  then $\pi(i+1i)\neq m$. (The subscript `a'  in $\mathcal{M}_{a}$ stands for alternating).
 
 We refer to the $\LL_{n}(\nu)$ as $n$-fold self-intersection local  times.
 If we now set
 \begin{equation}
\KK_{n}(\nu)= \sum_{\om\in\mathcal{I}_{\al}} \LL_{n}(\nu)(\om),\label{21.5}
 \end{equation} 
that is  we add together the $n$-fold self-intersection local  times $\LL_{n}(\nu)(\om)$ for each $\om\in\mathcal{I}_{\al}$, then using the moment formula for Poisson processes we find that the $\KK_{n}(\nu)$
  satisfy moment formulas similar to (\ref{rilt.15intro}) except that for each partition $B_{1}\cup \cdots\cup B_{j}=[1,n]$ and each $m$, there will be   some $l$ with $\pi^{ -1}( m)\in B_{l}$.
  
  We now construct an intersection local time involving intersections between different paths in $\mathcal{I}_{\al}$. We follow \cite[Section 7]{LMR2}, but the situation here is easier since all $\LL_{n}(x,r)$ are $\mu_{m}$ integrable.  For any set $A$, we  use $S_{n}(A)\subset A^{n}$ to denote the  subset of $A^{n}$ with distinct entries. That is, if   $(a_{i_{1}},\ldots, a_{ i_{n}})\in S_{n}(A)$ then $a_{i_{j}}\neq a_{i_{k}}$ for $i_{j}\neq i_{k}$. Let
   \begin{equation}
\KK_{l_{1},\ldots, l_{n}}(x,r)=  \sum_{ (\om_{i_{1} },\ldots,\om_{i_{n} })\in S_{n}(\mathcal{I}_{\al}) }\,\,\prod_{j=1}^{n}\LL_{l_{j}}(x,r)(\om_{i_{j} }).\label{21.6}
 \end{equation} 
  It follows as in the proof of \cite[Theorem 7.1]{LMR2} that if $l=l_{1}+\cdots+l_{n}$ and $ \nu\in \mathcal{G}_{K}^{2l}$, Then
\begin{equation}
\KK_{l_{1},\ldots, l_{n}}(\nu):=\lim_{r\rar 0}\int \KK_{l_{1},\ldots, l_{n}}(x,r)\,d\nu (x)\label{10.24}
\end{equation}
exists in all $L^{p}(   \PP_{ \al})$. Note that $\KK_{l_{1},\ldots, l_{n}}(\nu)$ involves $n$ distinct paths in 
  $\mathcal{I}_{\al}$. As in \cite[Theorem 7.2]{LMR2} we have 
  \bt\label{theo-multiloop} For $\nu\in\mathcal{G}_{K}^{2l}$  
\be
L_{n}(\nu)= \sum_{D_{1}\cup\cdots\cup D_{l}=[1,n]}
\KK_{|D_{1}|,\ldots, |D_{l}|}(\nu),\label{10.26}
\ee
where the sum is over all partitions of $[1,n]$. 
\et

\def\noopsort#1{} \def\printfirst#1#2{#1}
\def\singleletter#1{#1}
            \def\switchargs#1#2{#2#1}
\def\bibsameauth{\leavevmode\vrule height .1ex
            depth 0pt width 2.3em\relax\,}
\makeatletter
\renewcommand{\@biblabel}[1]{\hfill#1.}\makeatother
\newcommand{\bysame}{\leavevmode\hbox to3em{\hrulefill}\,}

 \def\wh{\widehat}
\def\ol{\overline}

\begin{thebibliography}{10}


\bibitem{Arc-Gin} M.~Arcones and E.~Gin\'e.
\newblock On decoupling, series expansion, and tail behavior of chaos
processes.
\newblock {\em Jour.\ Theoret.\ Prob.}, 6:101--122, 1993.

\bibitem{BKb} R.F. Bass and D. Khoshnevisan, Intersection local times and
Tanaka formulas. {\sl  Ann.  Inst. H. Poincar\'e \bf 29} (1993) 419--451.


\bibitem{BR} R.F. Bass and J.~Rosen,  An almost sure invariance principle
    for the range of  planar random walks, {\sl Ann. Probab.}, to appear.

\bibitem {DM4}
C. Dellacherie,  B. Maisonneuve and P.-A. Meyer,  (1992).
\newblock {\em Probabilities et Potential, Chapitres XVII a XXIV}.
\newblock Paris: Hermann.

\bibitem{Da} E. B.~Dynkin,  Regularized self-intersection local times of planar
Brownian motion. {\sl Ann. Probab. \bf 16} (1988) 58--74.

 \bibitem {five}
N. Eisenbaum, H. Kaspi, M. Marcus, J. Rosen and Zhan Shi, 
A Ray-Knight theorem for symmetric Markov processes,
 {\it Ann. Probab.},\,
 {\bf 28}\, (2000), 1781-1796.




\bibitem{LGa} J.-F.~Le Gall,  Wiener sausage and self intersection  local times.
{\sl J. Funct.\ Anal. \bf  88} (1990) 299--341.

\bibitem{LGb} J.-F. {Le Gall}.
        Propri\'et\'es d'intersection des marches al\'eatoires, I.
          Convergence vers le temps local d'intersection,
     {\sl Comm.\ Math.\ Phys.  \bf  104} (1986) 471--507.

 

%\bibitem{LMR1} Y. Le Jan, M. B.  Marcus and J.~Rosen,
%\newblock  Permanental fields,  loop soups and continuous additive functionals.
%\newblock  http://arxiv.org/pdf/1209.1804.pdf 

\bibitem{LMR2} Y. Le Jan, M. B.  Marcus and J.~Rosen,
\newblock  Intersection local times, loop soups and permanental Wick powers.
\newblock http://arxiv.org/pdf/1308.2701.pdf
 

\bibitem{K} J. F. C. Kingman, {\em Poisson Processes}, Oxford Studies in Probability, Clarendon Press, Oxford, (2002).

\bibitem{GR} I. Gradshteyn and I. Ryzhik, {\em Table of Integrals, Series and Products}, Academic Press, Oxford, (1980).

 

 \bibitem{MRmem} M. B.  Marcus and J.~Rosen,  Renormalized self-intersection local times and Wick power chaos processes, 
  \,\,{\it Memoirs of the A.M.S.},\, (1999), Number 675.

 
 \bibitem{MR96} M. B.  Marcus and J.~Rosen, Gaussian chaos and sample path properties of additive functionals of
  symmetric {M}arkov processes,
 {\it Ann. Probab.}\, {\bf 24}\, (1996), 1130--1177.
 
\bibitem{book} M. B.  Marcus and J.~Rosen, {\em Markov Processes,
Gaussian Processes and Local Times}, Cambridge University Press, New
York,  (2006).


   \bibitem{R} J.~Rosen, Joint continuity and a {D}oob-{M}eyer type decomposition for renormalized intersection
local times.{\it Ann. Inst. Henri Poincare },\,{\bf 35}\, (1999), 143--176.


%\bibitem{Ra} J.~Rosen,  Random walks and intersection local time.
     %   {\sl Ann. Probab. \bf  18} (1990) 959--977.


\bibitem{Rb} J.~Rosen,  Joint continuity of renormalized intersection  local
times, {\sl Ann.\
         Inst.\ H.~Poincar{\'e} \bf  32} (1996) 671--700.





\bibitem{Sz1} A.-S. Sznitman,
\newblock  Topics in occupation times and Gaussian free fields.
 {\it Zurich Lectures in Advanced Mathematics}, EMS, Zurich, 2012.
 
 \bibitem{Sz2} A.-S. Sznitman,
\newblock  An isomorphism theorem for random interlacements.
 {\it ECP},  {\bf 17} , 2012, no. 9, 1-9.
 
  \bibitem{Sz3} A.-S. Sznitman,
\newblock  On scaling limits and Brownian interlacements.

\bibitem{Va} S.~R.~S. Varadhan.
        Appendix to {E}uclidean quantum field theory by {K}. {S}y\-man\-zyk.
         In R.~Jost, editor,  {\sl Local Quantum Theory}. Academic  Press,
Reading, MA,
       1969.




\end{thebibliography}
\end{document}